\documentclass[journal]{IEEEtran}
\usepackage{amsmath,amsfonts,amsthm}
\usepackage[ruled,commentsnumbered,linesnumbered]{algorithm2e}
\usepackage{array}
\usepackage{textcomp}
\usepackage[dvipsnames]{xcolor}
\usepackage{stfloats}
\ifCLASSOPTIONcompsoc
\usepackage[caption=false,font=normalsize,labelfont=sf,textfont=sf]{subfig}
\else
\usepackage[caption=false,font=footnotesize]{subfig}
\fi
\usepackage{url}
\usepackage{verbatim}
\usepackage{graphicx}
\usepackage{cite}
\usepackage[colorlinks,linkcolor=BrickRed,citecolor=Aquamarine,urlcolor=cyan]{hyperref} 
\usepackage{makecell}
\usepackage{multirow}
\usepackage{color}
\usepackage{booktabs}
\usepackage{amsmath}
\usepackage{colortbl}
\usepackage{bbm}
\usepackage{textcomp}
\usepackage{xcolor}

\allowdisplaybreaks[4]
\usepackage{amssymb}
\usepackage{threeparttable}

\newtheorem{lemma}{Lemma}
\newtheorem{definition}{Definition}

\newtheorem{theorem}{Theorem}
\newtheorem{corollary}{Corollary}

\newtheorem{model}{Model}
\definecolor{lightblue}{rgb}{0.93,0.95,1.0}

\begin{document}

\title{Statistical Rounding Error Analysis for Random Matrix Computations}
\author{Yiming Fang,~\IEEEmembership{Graduate Student Member,~IEEE}, and
    Li Chen,~\IEEEmembership{Senior Member,~IEEE}
   
    \thanks{Yiming Fang and Li Chen are with the CAS Key Laboratory of Wireless-Optical Communications, University of Science and Technology of China, Hefei 230027, China (e-mail: fym1219@mail.ustc.edu.cn; chenli87@ustc.edu.cn).}
}



\maketitle
\begin{abstract}
The conventional rounding error analysis provides worst-case bounds with an associated failure probability and ignores the statistical property of the rounding errors. In this paper, we develop a new statistical rounding error analysis for random matrix computations. Such computations have numerous applications in the field of wireless communications, signal processing, and machine learning. By assuming the relative errors are independent random variables, we derive the approximate closed-form expressions for the expectation and variance of the rounding errors in various key computations for random matrices. Numerical experiments validate the accuracy of our derivations and demonstrate that our analytical expressions are generally at least two orders of magnitude tighter than alternative worst-case bounds, exemplified through the inner products. 
\end{abstract}

\begin{IEEEkeywords}
Floating-point arithmetic, IEEE 754 standard, Matrix computations, Random matrices, Rounding error analysis
\end{IEEEkeywords}

\section{Introduction}
\label{sec:intro}
\IEEEPARstart{R}{ounding} error analysis is a crucial method for assessing the numerical stability of algorithms, aiming to refine it based on the intrinsic properties obtained \cite{11076156,wilkinson1971modern,709374}. Classical rounding error analysis obtains backward error bounds involving the constant $\gamma_n=nu/\left(1-nu\right)$ for a dimension $n$ and unit roundoff $u$ \cite[Chapter 3]{higham2002accuracy}, which offer reasonable backward errors and valuable insights for double-precision arithmetic and moderate dimensions $n$. Nevertheless, the above classical rounding error analysis is very pessimistic about the problem with large dimensions and low-precision arithmetic \cite{higham2022mixed}.

To address the pessimistic problem of the classical rounding error analysis, probabilistic rounding error analysis has garnered significant attention. This method models relative errors $\delta_i$ as random variables to provide more accurate estimates of their average behavior. Neumann and Goldstine \cite{von_neumann} linearized the forward error $p$ as a sum, utilizing the central limit theorem to derive the probability distribution of $p$. To avoid the assumption that $n$ is sufficiently large for using the central limit theorem, Higham and Mary \cite{doi:10.1137/18M1226312} utilized a concentration inequality to yield the bounds proportional to $\sqrt{n}u$, but its bounds are still pessimistic. By assuming the data and the relative errors are random variables, Higham and Mary further derived tighter bounds and proposed a new algorithm for matrix-matrix products in \cite{doi:10.1137/20M1314355}. Moreover, Ipsen and Zhou \cite{doi:10.1137/19M1270434} derived forward error bounds for inner products with a clear relationship between failure probability and relative error. Additionally, the probabilistic forward error bounds for three classes of mono-precision summation algorithms were derived in \cite{hallman2023precision}.

Note that the aforementioned works conducted the rounding error analysis from a worst-case bound perspective. The authors in \cite{10845867} applied the above worst-case bounds to the communication precoding and detection processes and found them to be overly pessimistic. This is because the communication precoding and detection processes involve random matrix computations. In this case, worst-case scenarios are often rare from a statistical perspective. In other words, when the input is a random variable, both the output after rounding and the associated rounding error are also random variables, which is more reasonable to consider the statistical property of rounding errors rather than worst-case bounds. Therefore, if we can conduct the rounding error analysis from a statistical perspective, then we may obtain more accurate and tighter results for the rounding errors of random matrix computations. Random matrix computations have numerous applications in the field of wireless communications, signal processing, and machine learning.

 
In wireless communications, channel matrices are typically treated as random vectors or matrices \cite{tulino2004random}. Furthermore, precoding and detection processes based on the channel matrix involve random matrix-matrix computations and matrix factorization \cite{6457363,5730587,1261332}. In signal processing, random matrix computations are integral to covariance estimation algorithms \cite{5484583,10246801} and radar waveform design \cite{6649991,10547068}. In addition to wireless communications and signal processing, the computation of random matrices is also crucial in various machine learning tasks \cite{couillet2022random, NIPS2017_0f3d014e}.

The above applications motivate the study to give statistical rounding error analysis and derive the statistical property of rounding errors for the computation of random matrices. There is little research on how to obtain much sharper and more accurate results from a statistical perspective. Note that Constantinides et al. \cite{9048893} and Dahlqvist et al. \cite{constantinides2021rigorous} derived closed-form expressions for the distribution of the rounding errors under scalar computation when input data are random variables. However, to the best of our knowledge, the expectation and variance of the rounding errors for random matrix computations are still unknown. Compared to scalar computation, vector and matrix computation involve higher dimensions and different element combinations, leading to the progressive accumulation of rounding errors.

Motivated by the above observations, in this paper, we give a statistical rounding error analysis for random matrix computations and derive approximate closed-form expressions for the expectation and variance of the rounding errors. Our main contributions are summarized as follows:
\begin{itemize}
    \item \textbf{General rounding error analysis for random matrices.} We analyze rounding errors for the computation of random matrices with unknown specific distributions from a statistical perspective. Specifically, we derive the expectation and variance of the rounding errors for inner products and provide their approximate closed-form expressions by rigorous proof. Moreover, our analytical expressions can reduce to the probabilistic bounds in \cite{doi:10.1137/19M1270434,doi:10.1137/18M1226312} based on Bienaym\'e--Chebyshev inequality. Building on the results for inner products, we extend our analysis to rounding errors in matrix-vector and matrix-matrix products.
    \item \textbf{Specific rounding error analysis for Wishart matrices.} We conduct a statistical rounding error analysis for Wishart matrices by utilizing zero-forcing (ZF) detection and corresponding least squares (LS) problems as examples. More specifically, we present the rounding error analysis for standard algorithms used in solving LS problems, such as matrix factorization and the solution of triangular systems, and derive corresponding approximate closed-form expressions under the condition of Wishart matrices.
    \item \textbf{Tighter analytical expressions.} Most existing literature has almost exclusively focused on deriving the worst-case bounds, encompassing both deterministic loose bounds \cite{higham2002accuracy,doi:10.1137/19M1270434} and probabilistic bounds depending on a pessimistic failure probability to ascertain their validity \cite{doi:10.1137/18M1226312,doi:10.1137/20M1314355,doi:10.1137/19M1270434}. Moreover, some of these bounds are exact only in the first-order term and contain higher-order terms such as ``$+\mathcal{O}\left(u^2\right)$" \cite{doi:10.1137/20M1314355}. Many worst-case bounds do not yield closed-form expressions. In contrast, our analytical expressions are much tighter.
    To further demonstrate the superiority of our analytical expressions, the mean square error (MSE) of the rounding errors serves as the metric for comparing our analytical expressions with other worst-case bounds \cite{higham2002accuracy,doi:10.1137/19M1270434,doi:10.1137/18M1226312,doi:10.1137/20M1314355}. Numerical experiments demonstrate that our analytical expressions are generally at least two orders of magnitude tighter than alternative worst-case bounds, exemplified through the inner products. 
\end{itemize}

{Organization:} The paper is organized as follows. We first provide the probabilistic floating-point arithmetic model in Section \ref{sec:pro_fpm}. Then, we apply this model in Section \ref{sec:app_la} to present the general statistical rounding error for random matrices. Furthermore, a specific statistical rounding error for Wishart matrices is provided in Section \ref{sec:wishart}.
In Section \ref{sec:experiments}, we conduct a series of numerical experiments to validate our analytical expressions thoroughly. We summarize the conclusions in Section \ref{sec:conclusions}.

{Notation:} Bold uppercase letters denote matrices and bold lowercase letters denote vectors. For a matrix $\bf A$, ${\bf A}^T$, ${\bf A}^H$ and ${\bf A}^{-1}$ denote the transpose, the Hermitian transpose and inverse of ${\bf A}$, respectively. ${a}_{ij}$ denotes $(i,j)$th entry of ${\bf A}$. $\mathbb{E}\left(\cdot\right)$ and $\mathbb{V}\left(\cdot\right)$ denotes the expectation and variance, respectively. $ \left| {\bf A}\right|$ represents the matrix of absolute values $(\left| {a}_{ij}\right|)$. Given random variable $a$ and $b$, $\mathbb{C}(a,b)$ are the covariance between $a$ and $b$.

\section{Floating-point arithmetic model}
\label{sec:pro_fpm}
We first recall some basic definitions of floating-point arithmetic. A floating-point number system $\mathbb{F}$ is a subset of real numbers whose elements can be expressed as \cite{higham2002accuracy}
\begin{equation}
\label{eq:fp_representation}
    f = \pm m\times \eta^{e-t+1},
\end{equation}
where $\eta = 2$ is the base, the integer $t$ is the precision, the integer $e$ is the exponent within the range $e_{\min}\leq e\leq e_{\max}$, and the integer $m$ is the significand satisfying $0\leq m \leq \eta^t - 1$. Table \ref{tab:parafp} provides parameters for four floating-point arithmetic systems according to the IEEE standard \cite{4610935}.

\begin{table}[t]
\centering 
\setlength{\tabcolsep}{3.5pt}
\caption{Parameters for Four Floating-point Arithmetic}
\label{tab:parafp}
\begin{threeparttable}
\begin{tabular}{cclll}
\toprule[1pt]
\midrule
 & $\left(\mathrm{sig.}, \mathrm{exp.}\right)$\tnote{(1)} & \multicolumn{1}{c}{$u$\tnote{(2)}} & \multicolumn{1}{c}{$x_{\min}$\tnote{(3)}} & \multicolumn{1}{c}{$x_{\max}$\tnote{(4)}} \\ \midrule
$\mathtt{bfloat16}$    &  $(8,8)$    & $3.91\times 10^{-3}$ & $1.18\times 10^{-38}$& $3.39\times 10^{38}$ \\ 
$\mathtt{fp16}$      &  $(11,5)$    & $4.88\times 10^{-4}$ & $6.10\times 10^{-5}$& $6.55\times 10^{4}$ \\ 
$\mathtt{fp32}$     &  $(24,8)$    & $5.96\times 10^{-8}$ & $1.18\times 10^{-38}$& $3.40\times 10^{38}$ \\ 
$\mathtt{fp64}$     &  $(53,11)$    & $4.88\times 10^{-16}$ & $2.22\times 10^{-308}$& $1.80\times 10^{308}$ \\ \midrule
\bottomrule[1pt]
\end{tabular}
    \begin{tablenotes}    
        \footnotesize               
        \item[(1)] $\left(\mathrm{sig.}, \mathrm{exp.}\right)$ represents number of bits in significand and exponent.          
        \item[(2)] $u=\frac{1}{2}\eta^{1-t}$ is the unit roundoff.
        \item[(3)] $x_{\min}$ is the smallest normalized positive number.
        \item[(4)] $x_{\max}$ is the largest finite number.
    \end{tablenotes}            
\end{threeparttable}
\end{table} 

Then, considering round to nearest (RTN), the standard floating-point arithmetic can be given in the following model.
\begin{model}[Standard floating-point arithmetic model \cite{higham2002accuracy}]\label{model:det}
    Denote by $u$ the {unit roundoff}. The floating-point system $\mathbb{F}$ adheres to a standard arithmetic model if, for any $x, y \in \mathbb{F}$, one has 
    \begin{align}
    \label{eq:det_model}
	\boldsymbol{fl}\left( x\,\mathrm{op}\,y \right) =\left( x\,\mathrm{op}\,y \right) (1+\delta ),~\left| \delta \right|\le u,
    \end{align}
    where $\mathrm{op}\in \left\{ +,-,\times ,/ \right\}$, and $\boldsymbol{fl}\left( x\,\mathrm{op}\,y \right)$ is the correctly rounded (to nearest) value of $x\,\mathrm{op}\,y$.  
\end{model}
Model \ref{model:det} illustrates that finite-precision arithmetic operations can introduce a relative error $\delta$ for some $\delta \in [-u,u]$. In other words, Model \ref{model:det} is non-deterministic, and the specific value of $\delta$ in \eqref{eq:det_model} is unknown. This means that numerical analysis based on Model \ref{model:det} must consider all possible values $\delta$, i.e. it is based on Model \ref{model:det} is fundamentally a worst-case analysis. Notably, when the inputs are random variables (i.e., probabilistic), both the output after rounding and the associated relative error $\delta$ become random variables from a statistical perspective. In such cases, Model \ref{model:det} may be overly conservative, as worst-case scenarios are often rare and may not occur in practice. To this end, the following probabilistic model of relative error and probabilistic floating-point arithmetic is provided.


\begin{model}[Probabilistic model of relative error]\label{model:pro}
    Let the input signal be independent random variables. In the computation process, the relative errors $\delta$ in Model \ref{model:det} associated with every pair of operands are assumed to be independent random variables sampled from a given distribution $\mathcal{DIST}$ and satisfy $\left| \delta \right|\le u$, where the probability density function (PDF) of $\delta$ is given by {\rm \cite{9048893,constantinides2021rigorous}}
    \begin{equation}
        \label{eq:pdf_delta}
        f_{\delta}\left( t \right) \approx \begin{cases}
	\frac{3}{4u}&		t\in \left[ -\frac{u}{2},\frac{u}{2} \right]\\
	\frac{1}{2u}\left( \frac{u}{t}-1 \right) +\frac{1}{4u}\left( \frac{u}{t}-1 \right) ^2&	t\in \left[ -u,-\frac{u}{2} \right) \cup \left( \frac{u}{2},u          \right]\\
        \end{cases}.
    \end{equation}
\end{model}

Compared to \cite[Model 2]{doi:10.1137/20M1314355}, we do not need the input data to be bounded. It is interesting to note from {Model} \ref{model:pro} that the distribution of the relative error $\delta$ in \eqref{eq:pdf_delta} is approximately deterministic and is assumed to be independent of the input distribution. Furthermore, we can calculate the expectation and variance of $\delta$ using \eqref{eq:pdf_delta} as follows:
\begin{align}
    \mathbb{E} \left( \delta \right) &\approx \int_{-u}^{-\frac{u}{2}}{t\left[ \frac{1}{2u}\left( \frac{u}{x}-1 \right) +\frac{1}{4u}\left( \frac{u}{t}-1 \right) ^2 \right] dt} \nonumber\\
    &+\int_{-\frac{u}{2}}^{\frac{u}{2}}{\frac{3t}{4u}dt}+\int_{\frac{u}{2}}^u{\frac{t}{2u}\left( \frac{u}{t}-1 \right) dt} \nonumber\\
    &+\int_{\frac{u}{2}}^u{\frac{t}{4u}\left( \frac{u}{t}-1 \right) ^2dt} \nonumber\\
    &=0\label{eq:delta_e},\\
    \mathbb{V} \left( \delta \right) &\approx \mathbb{E} \left( \delta ^2 \right) -\left[ \mathbb{E} \left( \delta \right) \right] ^2=\mathbb{E} \left( \delta ^2 \right)\nonumber \\
    &=\int_{-u}^{-\frac{u}{2}}{t^2\left[ \frac{1}{2u}\left( \frac{u}{t}-1 \right) +\frac{1}{4u}\left( \frac{u}{t}-1 \right) ^2 \right] dt}\nonumber\\
    &+\int_{-\frac{u}{2}}^{\frac{u}{2}}{\frac{3t^2}{4u}dx}+\int_{\frac{u}{2}}^u{\frac{t^2}{2u}\left( \frac{u}{t}-1 \right) dt}\nonumber\\
    &+\int_{\frac{u}{2}}^u{\frac{t^2}{4u}\left( \frac{u}{t}-1 \right) ^2dt}\nonumber\\
    &=\frac{1}{6}u^2\triangleq \sigma ^2 \label{eq:delta_var}. 
\end{align}

Moreover, based on {Model} \ref{model:pro}, we can obtain the probabilistic floating-point arithmetic in the following definition.
\begin{definition}[Probabilistic floating-point arithmetic \cite{constantinides2021rigorous}]\label{lem:dis_delta}
    The floating-point system $\mathbb{F}$ adheres to a probabilistic arithmetic model if, for any random variables $x, y \in \mathbb{F}$, we have 
    \begin{align*}
    \boldsymbol{fl}\left( x\,\mathrm{op}\,y \right)& =\left( x\,\mathrm{op}\,y \right) (1+\delta ) \\
    &= \left( x\,\mathrm{op}\,y \right) + \Delta,~\delta \sim \mathcal{DIST},~\mathrm{op}\in \left\{ +,-,\times ,/ \right\}
    \end{align*}
    where $\Delta =\left( x\,\mathrm{op}\,y \right)\delta $ is the rounding error, and the PDF of $\delta$ is shown in \eqref{eq:pdf_delta}.
\end{definition}
Based on {Model} \ref{model:pro} and \eqref{eq:delta_e}, we observe that the rounding error $\Delta$ is approximately uncorrelated with the input because
\begin{align*}
    \mathbb{E} \left\{ \left( x\,\mathrm{op}\,y \right) \Delta \right\}& =\mathbb{E} \left\{ \left( x\,\mathrm{op}\,y \right) ^2\delta \right\} =\mathbb{E} \left\{ \left( x\,\mathrm{op}\,y \right) ^2 \right\} \mathbb{E} \left( \delta \right) \\
    &\approx 0 \approx\mathbb{E} \left\{ \left( x\,\mathrm{op}\,y \right) \right\} \mathbb{E} \left( \Delta \right).
\end{align*}

{Model} \ref{model:pro} and corresponding {Definition} \ref{lem:dis_delta} are not always realistic. For instance, real-world inputs may lack randomness. Additionally, repeated pairs of operands may result in the same $\delta$, i.e., they are dependent. Nevertheless, the question is whether these assumptions effectively model the actual rounding errors encountered in our computations (See similar comments of Hull and Swenson \cite{hull1966tests} and Kahan \cite{kahanimprobability}, as discussed in \cite{doi:10.1137/18M1226312,doi:10.1137/20M1314355}). We will show that the outcomes obtained under {Model} \ref{model:pro} closely approximate the actual results through numerical experiments in Section \ref{sec:experiments}. Furthermore, numerical experiments show that {Model} \ref{model:pro} is valid for independent random variables.

Note that many rounding error analyses depend on {Model} \ref{model:det}, and hence can potentially obtain the distribution (or expectation and variance) of their rounding errors if they can utilize {Model} \ref{model:pro}. Moreover, {Definition} \ref{lem:dis_delta} only addresses the distribution of rounding errors for scalar computations. When dealing with vector and matrix computations, the complexity and tediousness of calculations make it challenging to derive their rounding error expectation and variance. 

Overall, our goal in the following sections is to derive closed-form expressions of expectation and variance for the computation of random matrices based on {Model} \ref{model:pro} and to validate the accuracy of these expressions.

\section{General rounding error analysis for random matrices}
\label{sec:app_la}
In this section, we conduct a general rounding error analysis for random matrices whose specific distributions are unknown by leveraging {Model} \ref{model:pro}. First, the rounding error analysis for inner products for random vectors is given. Then we extend the rounding error analysis to encompass random matrix-vector and matrix-matrix products. 

Prior to discussing the rounding errors, we present several useful lemmas concerning expectation, variance, and some distributions that will be utilized in the subsequent sections.

\begin{lemma}[Expectation and variance of products of random variables {\cite[Section 2]{frishman1975arithmetic}}]\label{lem:e_var_product}
    Let $x$ and $y$ be independent random variables. Then the expectation and variance of their products $xy$ are given by
    \begin{align*}
        \mathbb{E}\left(xy\right) &= \mathbb{E}\left(x\right)\mathbb{E}\left(y\right),\\
        \mathbb{V}\left(  xy \right) &= \mathbb{V}\left(  x \right)\mathbb{V}\left(  y \right)+\mathbb{V}\left( y \right)\left[ \mathbb{E}\left(x\right)\right]^2 + \mathbb{V}\left( x \right)\left[ \mathbb{E}\left(y\right)\right]^2,
    \end{align*}
    Further, if and only if $\mathbb{E}\left(x\right)=\mathbb{E}\left(y\right)=0$, we have
    \begin{equation*}
        \mathbb{V}\left(  xy \right)=\mathbb{V}\left(  x\right)\mathbb{V}\left(  y \right).
    \end{equation*}
\end{lemma}

\begin{lemma}[Expectation and variance of random variables satisfying Wishart distribution]\label{lem:e_var_chi-square}
    Let random variable $x$ follow the chi-square distribution $\chi_{m}^2$, where $m>0$ represents the degree of freedom (DoF), then
    \begin{align*}
        \mathbb{E}\left(x\right) = m,~\mathbb{V}\left(x\right) = 2m.
    \end{align*}
\end{lemma}

\begin{lemma}[Expectation and variance of random variables satisfying student distribution]\label{lem:e_var_student}
    Let random variable $x$ follow the student's t-distribution $\mathcal{T} _{m}$ with $m$ DoF, then
    \begin{align*}
        \mathbb{E}\left(x\right) &= 0,~m>1\\
        \mathbb{V}\left(x\right) &= \frac{m}{m-2},~m>2
    \end{align*}
    where the expectation and variance of $x$ does not exist for $m=1$ and $m\leq2$, respectively.
\end{lemma}

\begin{lemma}[Expectation and variance of elements in random matrices satisfying Wishart distribution {\cite[Theroem 3.3.15]{gupta2018matrix}}]\label{lem:e_var_wishart}
    Let symmetric positive matrix ${\bf A}\in \mathbb{R}^{n\times n}$ follow the Wishart distribution $\mathcal{W}_n\left(m,\Sigma_n\right)$ with DoF $m$, then
    \begin{align*}
        \mathbb{E}\left(a_{ij}\right) & = m\sigma_{ij}, \\
        \mathbb{C}\left(a_{ij}, a_{kl}\right) & = m\left(\sigma_{ik}\sigma_{jl}+\sigma_{il}\sigma_{jk}\right),
    \end{align*}
    where $\Sigma = \left(\sigma_{ij}\right)$.
\end{lemma}

\begin{lemma}[Bienaym\'e--Chebyshev inequality \cite{heyde2012ij}]
\label{lem:BC}
Let $x$ be a random variable with finite expected value and finite nonzero variance. For any real number $\alpha>0$,
\begin{align*}
    \mathrm{Prob}\left( \left| x-\mathbb{E} \left( x \right) \right|\le \alpha \sqrt{\mathbb{V} \left( x \right)} \right) \ge 1-\frac{1}{\alpha ^2}.
\end{align*}
\end{lemma}

\subsection{Rounding Error Analysis for Inner Products}
We apply {Model} \ref{model:pro} to compute the inner product of two random vectors, and the expectation and variance of the rounding error are derived in the following theorem.
\begin{theorem}[Inner products]\label{the:inner}
    Let ${\bf x}, {\bf y} \in \mathbb{R}^{n\times1}$ be independent random vectors. The entries $x_i$ for $i=1$ to $n$ are sampled from a distribution with mean $\mu_x$ and variance $\sigma_x^2$, while the entries $y_i$ for $i=1$ to $n$ are sampled from a distribution with mean $\mu_y$ and variance $\sigma_y^2$. If $s={\bf x}^T{\bf y}$ is computed in floating-point arithmetic under {Model} \ref{model:pro}, the expectation and variance of the rounding error $\Delta s$ are given in \eqref{eq:inner_e} and \eqref{eq:inner_var} at the bottom of the page.     
    \begin{figure*}[hb] 
    \centering 
    \hrulefill
    \vspace*{2pt} 
    \begin{align}
    \mathbb{E} \left( \Delta s \right) &= 0,\label{eq:inner_e}\\
        \mathbb{V}\left( \Delta s \right) &\approx \tau \left[ \left( 1+\sigma ^2 \right) ^n+\frac{\left( 1+\sigma ^2 \right) ^2\left[ \left( 1+\sigma ^2 \right) ^{n-1}-1 \right]}{\sigma ^2}-n \right] \nonumber \\
        &+2\mu _{x}^{2}\mu _{y}^{2}\left[ \frac{\left( 1+\sigma ^2 \right) ^2\left[ \left( 1+\sigma ^2 \right) ^{n-1}-1 \right]}{\sigma ^4}-\frac{\left( n-1 \right) \left( 1+\sigma ^2 \right)}{\sigma ^2}-\frac{n\left( n-1 \right)}{2} \right]\nonumber\\
        &\triangleq \hbar\left(\mu_x,\sigma_x,\mu_y,\sigma_y,n,\sigma \right)\label{eq:inner_var}.
    \end{align}
    \end{figure*}
    Further, an asymptotic approximation of \eqref{eq:inner_var} can be express as
    \begin{align}\label{eq:v_sigma}
        \mathbb{V}\left( \Delta s \right) \approx \frac{\tau}{2}n^2\sigma ^2 + \frac{\mu _{x}^{2}\mu _{y}^{2}}{3}n^3\sigma ^2,
    \end{align}
    where $\sigma$ is defined in \eqref{eq:delta_var}, and $\tau =\sigma _{x}^{2}\sigma _{y}^{2}+\sigma _{x}^{2}\mu _{y}^{2}+\sigma _{y}^{2}\mu _{x}^{2}+\mu _{x}^{2}\mu _{y}^{2}$.
\end{theorem}
\begin{IEEEproof}
 The proof is available in Appendix \ref{app:inner}.
\end{IEEEproof}


{Theorem} \ref{the:inner} reveals that the variance of the rounding error is correlated with the variance, mean, and dimension of the input, and with the precision. Specifically, when the input distribution and precision are constant, the rounding error variance, i.e., MSE, grows asymptotically cubic and square with the input dimension for non-zero mean and zero mean variables, respectively. Thus, we can adjust the inputs to have zero means to obtain more accurate numerical results like the authors in \cite{doi:10.1137/20M1314355}. Furthermore, for fixed input distribution and dimension, the rounding error variance tends to zero as the precision increases. 

It is interesting to note that we can derive some probabilistic bounds based on expectation and variance in {Theorem} \ref{the:inner} in the following corollary.
\begin{corollary}[Probabilistic bounds for inner products]\label{co:pb}
Let ${\bf x}, {\bf y} \in \mathbb{R}^{n\times1}$ be independent random vectors. The entries $x_i$ for $i=1$ to $n$ are sampled from a distribution with mean $\mu_x$ and variance $\sigma_x^2$, while the entries $y_i$ for $i=1$ to $n$ are sampled from a distribution with mean $\mu_y$ and variance $\sigma_y^2$. If $s={\bf x}^T{\bf y}$ is computed in floating-point arithmetic under {Model} \ref{model:pro}, the computed $\hat{s}$ satisfies
\begin{equation}
    \left| \hat{s}-s \right|\lesssim \sqrt{\frac{\frac{\tau}{2}n^2\sigma ^2+\frac{\mu _{x}^{2}\mu _{y}^{2}}{3}n^3\sigma ^2}{\eta}},
\end{equation}
with probability at least $1-\eta$. Further, if $\mu_x\neq0$ and $\mu_y\neq0$, the backward error is bounded by
    \begin{equation}\label{eq:bwd_n0}
        \frac{\left| \hat{s}-s \right|}{\left| \mathbf{x} \right|^T\left| \mathbf{y} \right|}\lesssim \frac{1}{\sqrt{\eta}}\mathcal{O} \left( \sqrt{n}u \right),
    \end{equation}
    with probability at least $1-\eta$. 
    
    If $\left|x_i\right|,\left|y_i\right|$ has bounds and $\mu_x=0$ or $\mu_y = 0$, the backward error is bounded by
    \begin{equation}\label{eq:bwd_0}
        \frac{\left| \hat{s}-s \right|}{\left| \mathbf{x} \right|^T\left| \mathbf{y} \right|}\lesssim \frac{1}{\sqrt{\eta}}\mathcal{O} \left( u \right),
    \end{equation}
    with probability at least $1-\eta$.
\end{corollary}
\begin{IEEEproof}
 The proof is available in Appendix \ref{app:co}.
\end{IEEEproof}
Corollary \ref{co:pb} shows that we can obtain a probabilistic $\mathcal{O}(\sqrt{n}u)$ bound in \eqref{eq:bwd_n0} (similar to the bound in \cite{doi:10.1137/18M1226312}) from the derived result in {Theorem} \ref{the:inner} for random vectors. Moreover, when random vectors have a mean of zero, the backward error bound in \eqref{eq:bwd_0} is of order $\mathcal{O}(u)$ and does not grow $n$ (similar to the analysis in \cite{doi:10.1137/20M1314355}). Therefore, the derived results in Theorem \ref{the:inner} are general and can reduce to some classic bounds in \cite[Eq. (3.6)]{doi:10.1137/18M1226312} and \cite[Eq. (3.3)]{doi:10.1137/20M1314355}.




\subsection{Rounding Error Analysis for Matrix-Vector and Matrix-Matrix Products}
Building on the rounding error analysis for inner products, we can obtain the following theorems for matrix-vector and matrix-matrix products, respectively.
\begin{theorem}[Matrix-vector products]\label{the:mv}
    Let ${\bf A}\in \mathbb{R}^{m\times n}$ and ${\bf b} \in \mathbb{R}^{n\times1}$ be independent. Assume that the elements of $\bf A$ are sampled from a distribution with mean $\mu_a$ and variance $\sigma_a^2$, while the entries $\bf b$ are sampled from a distribution with mean $\mu_b$ and variance $\sigma_b^2$. If ${\bf y}={\bf A}{\bf b}$ is computed in floating-point arithmetic under {Model} \ref{model:pro}, the expectation and autocorrelation matrix of the rounding error $\Delta {\bf y}$ are given by
    \begin{align}
        \mathbb{E}\left( \Delta {\bf y} \right) &= \mathbf{0}_{m\times 1},\label{eq:mv_e}\\
        \mathbf{R}_{\Delta \mathbf{y}}& \approx \mathrm{diag}\left( \hbar ,\cdots ,\hbar \right)\label{eq:mv_var},
    \end{align}
    where $\hbar = \hbar\left(\mu_a,\sigma_a,\mu_b,\sigma_b,n,\sigma \right)$.
\end{theorem}
\begin{IEEEproof}
    The vector $\bf y$ is obtained by $m$ inner products, i.e., $y_i = {\bf a}_i^T{\bf x}$, where ${\bf a}$ is the $i$th row of $\bf A$. Therefore, using {Theorem} \ref{the:inner}, we have 
    \begin{equation*}
    \begin{matrix}
        \hat{y}_i=\boldsymbol{fl}\left( \mathbf{a}_{i}^{T}\mathbf{b} \right) =\mathbf{a}_{i}^{T}\mathbf{b}+\Delta \hat{y}_i, &i\in \left\{ 1,\cdots ,m \right\},
    \end{matrix}
    \end{equation*}
    where 
    \begin{align*}
        \mathbb{E} \left( \Delta \hat{y}_i \right) &= 0,\\
        \mathbb{V}\left( \Delta \hat{y}_i \right) &= \hbar\left(\mu_a,\sigma_a,\mu_b,\sigma_b,n,\sigma \right),\\
        \mathbb{E} \left( \Delta \hat{y}_i \Delta \hat{y}_j \right) &= 0,\,\,\,\, i\neq j.
    \end{align*}
    Note that $\Delta {\bf y}=\left[\Delta \hat{y}_1,\Delta \hat{y}_2,\cdots,\Delta \hat{y}_n\right]$, and we can know that 
    \begin{align*}
        \mathbb{E}\left( \Delta {\bf y} \right) &= \mathbf{0}_{m\times 1},\\
        \mathbf{R}_{\Delta \mathbf{y}}& =\mathbb{E} \left\{ \Delta \mathbf{y}\Delta \mathbf{y}^T \right\} \approx \mathrm{diag}\left( \hbar ,\cdots ,\hbar \right).
    \end{align*}
    Therefore, {Theorem} \ref{the:mv} holds.
\end{IEEEproof}

\begin{theorem}[Matrix-matrix products]\label{the:mm}
    Let ${\bf A}\in \mathbb{R}^{m\times n}$ and ${\bf B} \in \mathbb{R}^{n\times p}$ be independent. Assume that the elements of $\bf A$ are sampled from a distribution with mean $\mu_a$ and variance $\sigma_a^2$, while the entries $\bf B$ are sampled from a distribution with mean $\mu_b$ and variance $\sigma_b^2$. If ${\bf C}={\bf A}{\bf B}$ is computed in floating-point arithmetic under {Model} \ref{model:pro}, the expectation and autocorrelation matrix of the rounding error $\Delta {\bf C}$ are given by
    \begin{align}
        \mathbb{E}\left( \Delta {\bf C} \right) &= {\bf 0}_{m\times p}\label{eq:mm_e},\\
        \mathbf{R}_{\Delta \mathbf{C}}& = \mathrm{diag}\left( p\hbar, \cdots, p\hbar \right)\label{eq:mm_var}.
    \end{align}
\end{theorem}
\begin{IEEEproof}
    The matrix $\bf C$ can be obtained by $p$ matrix-vector products, i.e., ${\bf c}_j = {\bf A}{\bf b}_j$, where ${\bf c}_j$ and ${\bf b}_j$ are the $j$th row of $\bf C$ and $\bf B$, respectively. Then the following proof is similar to that of {Theorem} \ref{the:mv}, omitted for conciseness.
\end{IEEEproof}

Similar to {Theorem} \ref{the:inner}, {Theorem} \ref{the:mv} and \ref{the:mm} also shows the influence of four crucial factors on the variance of the rounding error: the mean, variance, and dimension of the input, and the precision.


\section{Specific rounding error analysis for Wishart Matrices}
\label{sec:wishart}
In this section, we provide the specific rounding error analysis for Wishart matrices. Wishart matrices arise when a Gram matrix is generated from a matrix with a Gaussian distribution, which has been extensively studied and applied to wireless communications \cite{1237156,1192179,5474635}. 

In the following presentation, we use ZF detection and the corresponding LS problem as an example and give the statistical rounding error analysis for the standard algorithms involved in the process of solving the problem.

Considering ZF detection in the uplink, the received vector after using the ZF detector $\bf x$ is given by
\begin{equation}
\label{eq:zf_c}
    {\bf x} = \left ({\bf H}^T{\bf H} \right)^{-1}{\bf H}^T{\bf z},
\end{equation}
where ${\bf H}\in \mathbb{R}^{m\times n}$ is the channel matrix with the entries $h_{ij}\sim \mathcal{N}(0,1)$ and $\bf z$ is the received signal at the base station. To avoid matrix inversion, we transform \eqref{eq:zf_c} into LS problem, yielding
\begin{equation}
\label{eq:ne}
    {\bf H}^T{\bf H} {\bf x} = {\bf H}^T{\bf z}. 
\end{equation}
One traditional approach for solving the LS problem is the normal equation method \cite[Alg. 5.3.1]{gloub1996matrix}. Specifically, we can use the following procedure:
\begin{equation}\label{eq:nep}
    \begin{aligned}
    &{\bf A} = {\bf H}^T{\bf H},~ {\bf c} = {\bf H}^T{\bf z},\\
    &{\bf A}  = {\bf R}^T{\bf R}~{\rm or}~{\bf A}  = {\bf L}{\bf U},\\
    &{\rm Solve}~{\bf R}^T{\bf y}={\bf c},~{\bf R}{\bf x}={\bf y}~{\rm or}~{\bf L}{\bf y}={\bf c},~{\bf U}{\bf x}={\bf y}.
    \end{aligned}
\end{equation}
where ${\bf A}$ is a Wishart matrix satisfying the Wishart distribution $\mathcal{W}_n\left(m,{\bf I}_n\right)$. It is observed that the standard algorithms such as matrix-matrix products, matrix-vector products, Cholesky factorization, LU factorization, and the solution of triangular systems are used in \eqref{eq:nep}. Since the rounding error analysis for matrix-matrix and matrix-vector multiplication has been presented in Section \ref{sec:app_la}, we further derive the rounding error analysis for the solution of triangular systems and matrix factorization under the condition of Wishart matrices.

\subsection{Rounding Error Analysis for The Solution of Triangular Systems} 
In the subsection, we provide the rounding error analysis for the solution of triangular systems. Without loss of generality, a lower triangular matrix is considered. Specifically, let the lower triangular matrix ${\bf T} = \left(t_{ij}\right)$ be obtained through the Cholesky factorization of Wishart matrix ${\bf A}\sim \mathcal{W}_n\left(m,{\bf I}_n\right)$ in \eqref{eq:nep}. Furthermore, we have $t_{ij}, 1\leq j\leq i\leq n$ are independently distributed, $t_{ii}^2\sim \chi_{m-i+1}^2,1\leq i\leq n$ and $t_{ij}\sim \mathcal{N}\left(0,1\right),1\leq j<i\leq n$ \cite[Theorem 3.3.4]{gupta2018matrix}. 


Then, given a triangular system ${\bf Tx}={\bf b}$, the solution $\bf x$ can be computed as follows: \cite[Algorithm 3.1.3]{GoVa13}
\begin{equation}
    \label{eq:tri_sub}
        x_i = \frac{b_i - \sum_{j=1}^{i-1}t_{ij}x_j}{t_{ii}}=\frac{b_i}{t_{ii}} -\sum_{j=1}^{i-1}\frac{t_{ij}x_j}{t_{ii}},~i=1:n.
\end{equation}
Using \eqref{eq:tri_sub} and {Model} \ref{model:pro}, we have the following theorem to obtain the expectation and variance of the rounding errors for the solution of triangular systems.

\begin{theorem}[Solution of triangular systems]\label{the:tri}
    Let ${\bf T} = \left(t_{ij}\right)\in \mathbb{R}^{n\times n}$ be a nonsingular and lower triangular matrix, where $t_{ij}, 1\leq j\leq i\leq n$ are independently distributed, $t_{ii}^2\sim \chi_{m-i+1}^2,1\leq i\leq n$ and $t_{ij}\sim \mathcal{N}\left(0,1\right),1\leq j<i\leq n$. Let the elements of ${\bf b}= \left(b_{i}\right)\in \mathbb{R}^{n\times 1}$ follow independent normal distribution $\mathcal{N}\left(0,1\right)$. If the triangular system ${\bf Tx}={\bf b}$ is solved by substitution in floating-point arithmetic under {Model} \ref{model:pro}, provided that $m>n+1$, the expectation and variance of the rounding error $\Delta x_i,1\leq i\leq n$ can be expressed as
    \begin{align}
        \mathbb{E}\left( \Delta x_i \right) &= 0\label{eq:tri_e},\\
        \mathbb{V}\left( \Delta x_i \right) &\approx\frac{\sum_{j=1}^{i-1}{\mathbb{V}\left( x_j \right) \left( 1+\sigma _{\psi _j}^{2} \right) \left( 1+\sigma ^2 \right) ^{i-j+2}}}{m-i-1}\nonumber\\
        &+\frac{\left( 1+\sigma ^2 \right) ^i}{m-i-1}-\mathbb{V}\left( x_i \right) \label{eq:tri_var},
    \end{align}
    where
    \begin{align}
        \sigma _{\psi _j}^{2}&=\frac{\mathbb{V}\left( \Delta x_j \right)}{\mathbb{V}\left( x_j \right)},~j=1:i-1,\\
        \mathbb{V}\left( x_i \right) &=\frac{1}{m-i-1}+\frac{1}{m-i-1}\sum_{j=1}^{i-1}{\mathbb{V}\left( x_j \right)}.
    \end{align}
\end{theorem}
\begin{IEEEproof}
    The proof is available in Appendix \ref{app:tri}.
\end{IEEEproof}

{Theorem} \ref{the:tri} indicates that the variance of the rounding error correlates with the previously derived outcomes. This correlation is reasonable, as solving triangular systems involves iterative processes, as depicted in \eqref{eq:tri_sub}. Additionally, when other parameters remain constant, the variance of the rounding error decreases as the DoF increases. This decrease occurs because a higher DoF results in an increased variance of $t_{ii}$.


\subsection{Rounding Error Analysis for LU Factorization}
In the subsection, we present the rounding error analysis for LU factorization of Wishart matrix ${\bf A}\sim \mathcal{W}_n\left(m,{\bf I}_n\right)$. Specifically, the LU factors of ${\bf A}$ are given by the Doolittle form of Gaussian elimination in the following recurrences \cite[Algorithm 9.2]{higham2002accuracy}:
\begin{equation}
    \begin{aligned}
\label{eq:d_lu}
             u_{kj}&=a_{kj}-\sum_{i=1}^{k-1}{l_{ki}u_{ij}},~j=k:n,\\
        l_{ik}&={{\left( a_{ik}-\sum_{j=1}^{k-1}{l_{ij}u_{jk}} \right)}\Bigg/{u_{kk}}},~i=k+1:n,
\end{aligned}
\end{equation}
where $1 \leq k\leq n$.
    
Before giving the expectation and variance of the lower triangular matrix ${\bf L}$ and the upper triangular matrix ${\bf U}$, we need to make some preparations. Specifically, the distribution of ${\bf L}\in \mathbb{R}^{n\times n}$, ${\bf U}\in \mathbb{R}^{n\times n}$, and the expected values of other terms in \eqref{eq:d_lu} are derived in the following two lemmas.

\begin{lemma}[Distribution of $\bf L$ and $\bf U$]\label{lem:dis_l_u}
    Given symmetric positive matrix ${\bf A} \sim \mathcal{W}_n\left(m,{\bf I}_n\right)$ and provided that $m>n+1$, the distribution of the LU factors of $\bf A$ can be summarized as follows:
    \begin{itemize}
        \item For the upper triangular matrix ${\bf U}=\left( u_{ij}\right)$, if $i=j$, we have
        \begin{equation*}
            u_{ii} \sim \chi_{\nu}^2,~ 1\leq i \leq n.           
        \end{equation*}
        If $i\neq j$, we have its PDF, i.e.,
        \begin{equation*}
            f_{u_{ij}}\left( z \right) =\begin{cases}
	   \frac{\Gamma \left( \frac{\nu -1}{2} \right)}{2\sqrt{\pi}\Gamma \left( \frac{\nu}{2} \right)},&		z=0\\
	   \frac{1}{\sqrt{2\pi}2^{\frac{\nu}{2}-1}\Gamma \left( \frac{\nu}{2} \right)}\left( \left| z \right| \right) ^{\frac{\nu -1}{2}}K_{\frac{\nu -1}{2}}\left(           \left| z \right| \right) ,&		else\\
            \end{cases},
        \end{equation*}
        where $ 1\leq i <j\leq n$, $\nu = m-i+1$, and $K_n\left(y\right)$ is the modified Bessel function of the second kind.
        \item For the lower triangular matrix ${\bf L}=\left( l_{ij}\right)$, we have
        \begin{equation*}
            l_{ij} = \frac{1}{\sqrt{m-j+1}}t_j,~ t_j\sim \mathcal{T}_{m-j+1},~ 1\le j<i\le n,
        \end{equation*}
    \end{itemize}
\end{lemma}
\begin{IEEEproof}
    For the Cholesky factorization ${\bf A} = {\bf R}^T{\bf R}$, where ${\bf R} = \left(r_{ij}\right)$ is a upper triangular matrix with $r_{ii}>0$, we have $r_{ij}, 1\leq i\leq j\leq n$ are independently distributed, $r_{ii}^2\sim \chi_{m-i+1}^2,1\leq i\leq n$ and $r_{ij}\sim \mathcal{N}\left(0,1\right),1\leq i<j\leq n$ \cite[Theorem 3.3.4]{gupta2018matrix}. Note that the relationship between LU factorization and Cholesky factorization can be established through the $LDL^T$ factorization, i.e.,
    \begin{equation*}
        \mathbf{U}=\mathbf{DL}^T=\mathbf{D}^{\frac{1}{2}}\mathbf{R},~\mathbf{L}=\mathbf{R}^T\mathbf{D}^{-\frac{1}{2}},~  \mathbf{D}=\mathrm{diag}\left( r_{11}^{2},\cdots ,r_{nn}^{2} \right).
    \end{equation*}

    Therefore, for the upper triangular matrix ${\bf U}=\left( u_{ij}\right)$, if $i=j$, we have
    \begin{equation*}
        u_{ii}=r_{ii}^{2}\sim \chi _{m-i+1}^{2},1\le i\le n.
    \end{equation*}
    Then, if $i\neq j$, we have $u_{ij}=r_{ii}r_{ij}$. Given $r_{ii}^2\sim \chi_{m-i+1}^2$ and $r_{ij}\sim \mathcal{N}\left(0,1\right)$, the joint PDF of 
    $r_{ii}$ and $r_{ij}$ can be expressed as
    \begin{align*}
        f_{r_{ii},r_{ij}}\left( x,y \right) &=f_{r_{ii}}\left( x \right) f_{r_{ij}}\left( y \right) \\
        &=\begin{cases}
	\frac{x^{\nu -1}e^{-\frac{x^2+y^2}{2}}}{\sqrt{2\pi}2^{\frac{\nu}{2}-1}\Gamma \left( \frac{\nu}{2} \right)},&		x>0\\
	0,&		x\le 0\\
        \end{cases},
    \end{align*}
    where $\nu = m-i+1$. Then the PDF of $u_{ij}$ is given by
    \begin{align*}
        f_{u_{ij}}\left( z \right) &=\int_{-\infty}^{+\infty}{\frac{1}{\left| x \right|}}f\left( x,\frac{z}{x} \right) dx=\int_0^{+\infty}{\frac{1}{x}}f\left( x,\frac{z}{x} \right) dx\\
        &\overset{\left( a \right)}{=}\frac{1}{\sqrt{2\pi}2^{\frac{\nu}{2}-1}\Gamma \left( \frac{\nu}{2} \right)}\left( \left| z \right| \right) ^{\frac{\nu -1}{2}}K_{\frac{\nu -1}{2}}\left( \left| z \right| \right), 
    \end{align*}
    where $(a)$ follows \cite[Eq. (3.478.4)]{gradshteyn2014table}, and $K_n\left(y\right)$ is the modified Bessel function of the second kind. Moreover, note that $K_n\left(y\right)\longrightarrow \infty$ when $y \longrightarrow 0$. We can use the asymptotic form for small arguments of the Bessel function \cite{abramowitz1948handbook} and have
    \begin{align*}
        f_{u_{ij}}\left( z \right) &=\frac{\left( \left| z \right| \right) ^{\frac{\nu -1}{2}}}{\sqrt{2\pi}2^{\frac{\nu}{2}-1}\Gamma \left( \frac{\nu}{2} \right)}\frac{\Gamma \left( \frac{\nu -1}{2} \right)}{2}\left( \frac{2}{\left| z \right|} \right) ^{\frac{\nu -1}{2}}\\
        &=\frac{\Gamma \left( \frac{\nu -1}{2} \right)}{2\sqrt{\pi}\Gamma \left( \frac{\nu}{2} \right)},~\left| z \right| \longrightarrow 0.
    \end{align*}
    
    For the lower triangular matrix ${\bf L}=\left( l_{ij}\right)$, we have
    \begin{align*}
        l_{ij}=\frac{r_{ji}}{r_{jj}}&=\frac{1}{\sqrt{m-j+1}}\left( \frac{\sqrt{m-j+1}r_{ji}}{r_{jj}} \right)\\
        &=\frac{1}{\sqrt{m-j+1}}t_j,~1\le j<i\le n.
    \end{align*}
    where $t_j\sim \mathcal{T}_{m-j+1}$.
    Hence, the lemma holds.
\end{IEEEproof}

{Lemma} \ref{lem:dis_l_u} reflects the fact that the off-diagonal elements in each row of $\mathbf{U}$ also share the same distribution, while the off-diagonal elements in each column of $\mathbf{L}$ share the same distribution.

    \begin{figure*}[hb] 
    \centering 
    \hrulefill
    \vspace*{2pt} 
    \begin{align}\small
        \mathbb{E}\left( \Delta u_{kk} \right) &= 0,\label{eq:e_ukk}\\
        \mathbb{V}\left( \Delta u_{kk} \right)&\approx \left( m^2-4 \right) \left[ \left( 1+\sigma ^2 \right) ^{k-1}-1 \right] -3\left( k-1 \right)+3\sum_{i=1}^{k-1}{\left( 1+\sigma _{\epsilon _i}^{2} \right) \left( 1+\sigma _{\eta _i}^{2} \right) \left( 1+\sigma ^2 \right) ^{k-i+1}} \nonumber\\
        &-2\left( m+2 \right) \left[ \frac{\left( 1+\sigma ^2 \right) \left[ \left( 1+\sigma ^2 \right) ^{k-2}-1 \right]}{\sigma ^2}-k+2 \right]. \label{eq:v_ukk} \\
            \mathbb{E}\left( \Delta u_{kj} \right) &= 0,\label{eq:e_ukj}\\
            \mathbb{V}\left( \Delta u_{kj} \right)&\approx\left( m-2 \right) \left[ \left( 1+\sigma ^2 \right) ^{k-1}-1 \right] +\sum_{i=1}^{k-1}{\left( 1+\sigma _{\epsilon _i}^{2} \right) \left( 1+\sigma _{\eta _i}^{2} \right) \left( 1+\sigma ^2 \right) ^{k-i+1}}
            -2\frac{\left( 1+\sigma ^2 \right) \left[ \left( 1+\sigma ^2 \right) ^{k-2}-1 \right]}{\sigma ^2}+k-3,\label{eq:v_ukj}\\
            \mathbb{E}\left( \Delta l_{ik} \right) &= 0,\label{eq:e_lik}\\
            \mathbb{V}\left( \Delta l_{ik} \right)&\approx \frac{(m-6)\left[ \left( 1+\sigma _{\eta _k}^{2} \right) \left( 1+\sigma ^2 \right) ^k-1 \right]}{\left( m-k-1 \right) \left( m-k-3 \right)}
            +\frac{\left( 1+\sigma _{\eta _k}^{2} \right) \sum_{j=1}^{k-1}{\left( 1+\sigma _{\epsilon _j}^{2} \right) \left( 1+\sigma _{\eta _j}^{2} \right) \left( 1+\sigma ^2 \right) ^{k-j+2}}-k+1}{\left( m-k-1 \right) \left( m-k-3 \right)} \nonumber \\
            &-\frac{2\left[ \left( 1+\sigma _{\eta _k}^{2} \right) \frac{\left( 1+\sigma ^2 \right) ^2\left[ \left( 1+\sigma ^2 \right) ^{k-2}-1 \right]}{\sigma ^2}-k+2 \right]}{\left( m-k-1 \right) \left( m-k-3 \right)}.\label{eq:v_lik}
    \end{align}
    \end{figure*}
    
\begin{lemma}[Expected values of other terms in \eqref{eq:d_lu}]\label{lem:term_lu_e_var}
    Given a symmetric positive matrix ${\bf A} \sim \mathcal{W}_n\left(m,{\bf I}_n\right)$ and provided that $m>n+3$, if the LU factorization of ${\bf A}$ is computed via \eqref{eq:d_lu}, we can have the expected values in \eqref{eq:d_lu} as follows:
    \begin{itemize}
        \item Denote $q_{kij} = l_{ki}u_{ij}, 1\leq i\leq k-1, k\leq j\leq n, 1\leq k\leq n$. If $j = k$, we have 
        \begin{align*}
            \mathbb{E} \left( q_{kij} \right) &=1,~ \mathbb{V}\left( q_{kij} \right) =2,~
            \mathbb{C}\left( q_{ksj},q_{ktj} \right) =0,~ s\neq t, \\
            \mathbb{C}\left( a_{kk},q_{kij} \right)& =\mathbb{E} \left( a_{kk}q_{kij} \right) =2.
        \end{align*}
        If $j\neq k$, we have
        \begin{align*}
            \mathbb{E} \left( q_{kij} \right)& =0,~ \mathbb{V}\left( q_{kij} \right) =1,~
            \mathbb{C}\left( q_{ksj},q_{ktj} \right) = 0,~ s\neq t, \\
            \mathbb{C}\left( a_{kj},q_{kij} \right)& =\mathbb{E} \left( a_{kj}q_{kij} \right)=1.
        \end{align*}
        \item Denote $p_{ijk}=\frac{l_{ij}u_{jk}}{u_{kk}}, k+1\leq i\leq n, 1\leq j\leq k-1, 1\leq k\leq n$. We have
        \begin{align*}
            \mathbb{E} \left( p_{ijk} \right) &=0,~ \mathbb{V}\left( p_{ijk} \right) =\frac{1}{\left( m-k-1 \right) \left( m-k-3 \right)},\\
            \mathbb{C}\left( p_{isk},p_{itk} \right) &=0.
        \end{align*}
        \item Denote $o_{ik}=\frac{a_{ik}}{u_{kk}}, k+1\leq i\leq n, 1\leq j\leq k-1, 1\leq k\leq n$. We have
        \begin{align*}
            \mathbb{E} \left( o_{ik} \right) &=0,~ \mathbb{V}\left( o_{ik} \right) =\frac{m-4}{\left( m-k-1 \right) \left( m-k-3 \right)},\\
            \mathbb{C}\left( o_{ik},p_{ijk} \right)&=\mathbb{E} \left( o_{ik}p_{ijk} \right) =\frac{1}{\left( m-k-1 \right) \left( m-k-3 \right)}. 
        \end{align*}
    \end{itemize}
\end{lemma}
\begin{IEEEproof}
    From the proof of {Lemma} \ref{lem:dis_l_u}, we derive $l_{ij}=\frac{r_{ji}}{r_{jj}}$ and $u_{ij}=r_{ii}r_{ij}$. Subsequently, we can obtain $q_{kij}=r_{ik}r_{ij}$, $p_{ijk}=\frac{r_{ji}r_{jk}}{r_{kk}^{2}}$, and $o_{ik}=\frac{r_{ki}+\sum_{j=1}^{k-1}{r_{ji}r_{jk}}}{r_{kk}}$. Note that the distribution of $r_{ij}$ is already known, thus concluding the proof.
\end{IEEEproof}

Next, based on the two lemmas above, the expectation and variance of the rounding errors for LU factors are presented in the following theorem.
\begin{theorem}[LU factorization]\label{the:lu}
    Given a symmetric positive matrix ${\bf A} \sim \mathcal{W}_n\left(m,{\bf I}_n\right)$ and provided that $m>n+3$, if the LU factorization of ${\bf A}$ is computed via \eqref{eq:d_lu} in float-point arithmetic under {Model} \ref{model:pro}, the expectation and variance of the rounding errors for LU factors can be expressed as follows:
    \begin{itemize}
        \item For the upper triangular matrix ${\bf U}=\left( u_{kj}\right)$, if $j=k$, the expectation and variance of the rounding errors $\Delta u_{kk}$ are given in \eqref{eq:e_ukk} and \eqref{eq:v_ukk} at the bottom of the page, respectively.
        
        If $j\neq k$, the expectation and variance of the rounding errors $\Delta u_{kj}$ are given in \eqref{eq:e_ukj} and \eqref{eq:v_ukj} at the bottom of the page, where $k\leq j \leq n$, $\sigma _{\epsilon _i}^{2}=\left( m-i-1 \right) \mathbb{V}\left( \Delta l_{ki} \right), \sigma _{\eta _i}^{2}=\frac{\mathbb{V}\left( \Delta u_{ij} \right)}{m-i+1}, 1\leq i\leq k-1,  1 \leq k\leq n$.
        \item For the lower triangular matrix ${\bf L}=\left( l_{ik}\right)$, the expectation and variance of the rounding errors $\Delta l_{ik}$ are given in \eqref{eq:e_lik} and \eqref{eq:v_lik} at the bottom of the page, where $k+1\leq i \leq n$, $\sigma _{\eta _k}^{2}=\frac{\mathbb{V}\left( \Delta u_{kk} \right)}{\left( m-k+1 \right) \left( m-k+3 \right)}, 1 \leq k\leq n$.
    \end{itemize}
\end{theorem}
\begin{IEEEproof}
    The proof is available in the Appendix \ref{app:lu}.
\end{IEEEproof}

{Theorem} \ref{the:lu} reveals that when other parameters are fixed, the variance of the rounding error for ${u}_{kj}$ increases with the DoF, whereas for ${l}_{ik}$, it decreases. This is because of the different distribution of the LU factors (see {Lemma} \ref{lem:dis_l_u}). Moreover, we can find that the variance of the rounding error is independent of the dimension $n$ provided that $m>n+3$.

\section{Simulation Results and Discussion}
\label{sec:experiments}
\begin{figure*}[t]
    \centering
    \subfloat[Random uniform $\mathcal{U}\left(0,1\right)$ vectors.]{\includegraphics[width=0.35\textwidth]{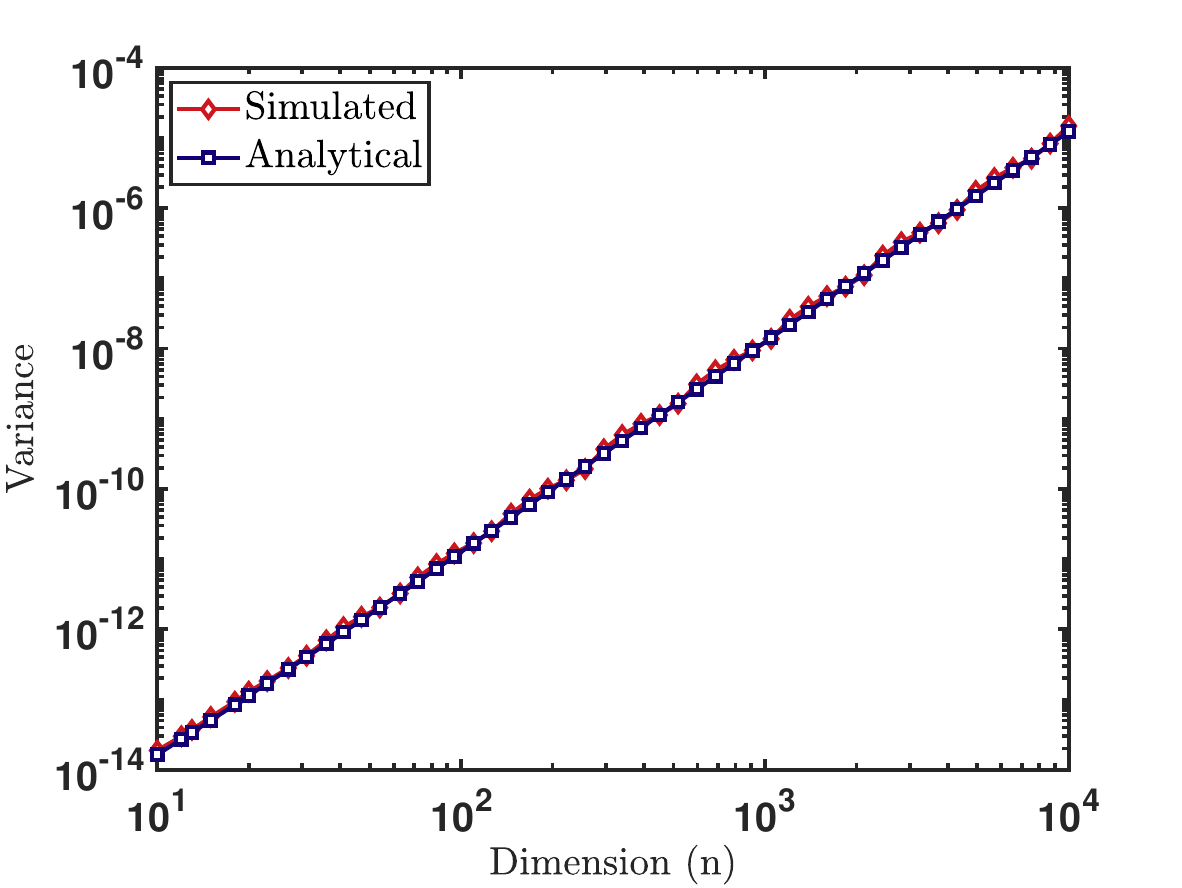}
    \label{fig:Inner_Product_Single_Uni_0_1}}
    \hspace{30pt}
    \subfloat[Random uniform $\mathcal{U}\left(-1,1\right)$ vectors.]{\includegraphics[width=0.35\textwidth]{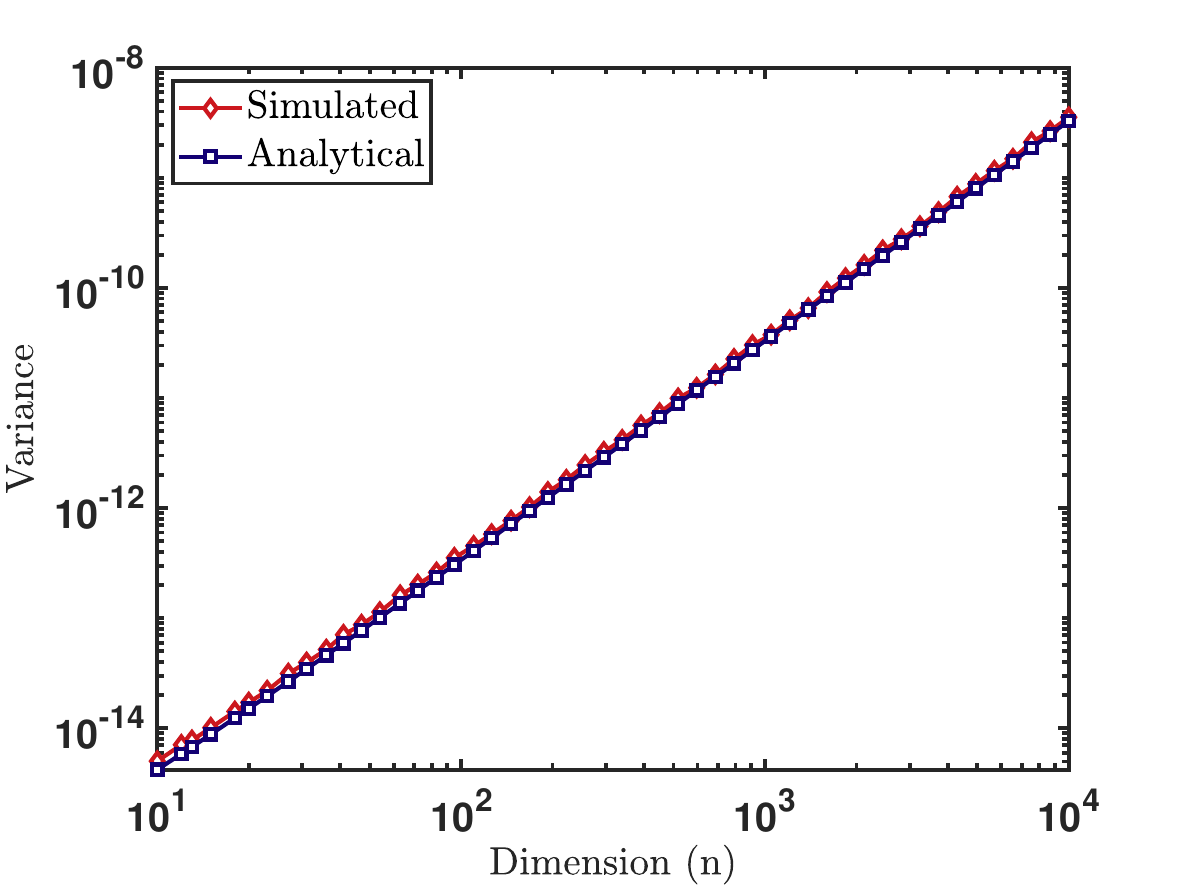}
    \label{fig:Inner_Product_Single_Uni_-1_1}}\\
    \subfloat[Random Gaussian $\mathcal{N}\left(0,1\right)$ vectors.]{\includegraphics[width=0.35\textwidth]{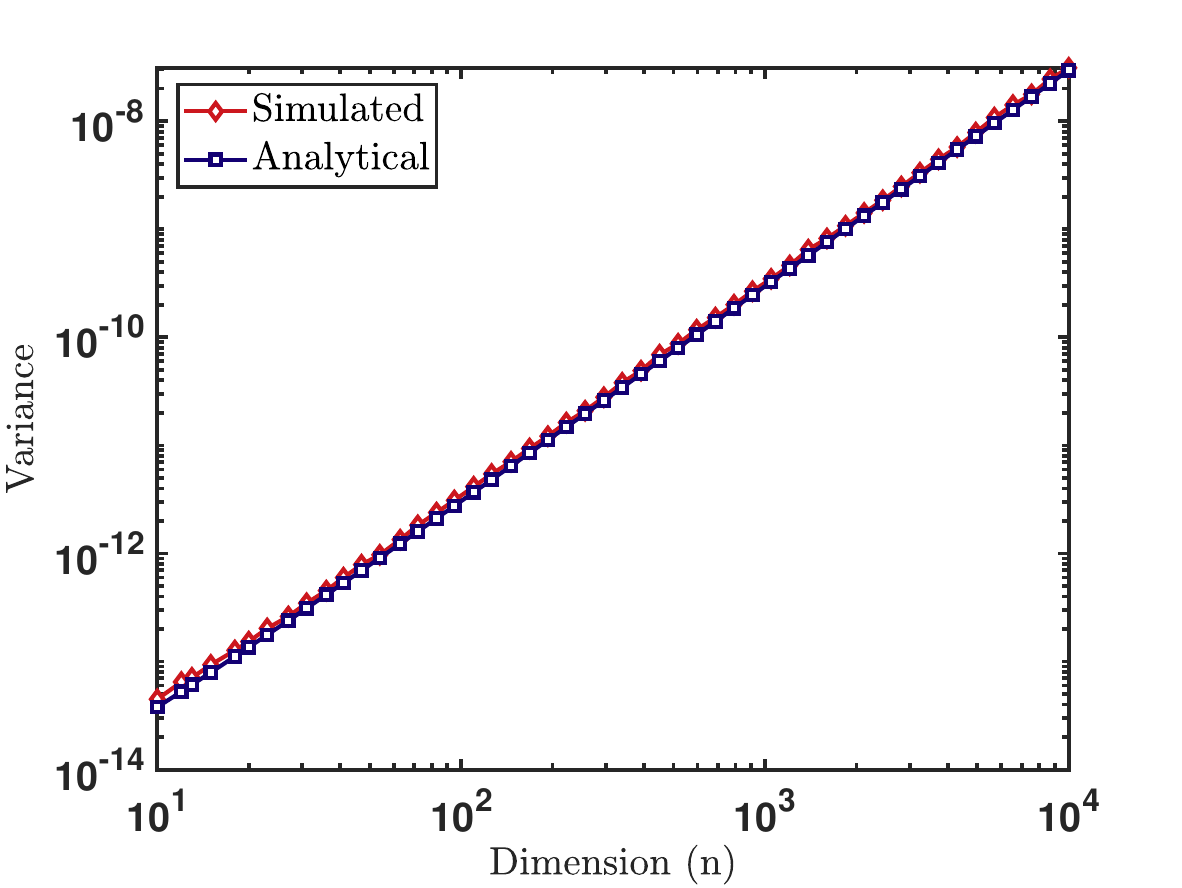}
    \label{fig:Inner_Product_Single_Nor_0_1}}
    \hspace{30pt}
    \subfloat[Random Gaussian $\mathcal{N}\left(1,1\right)$ vectors.]{\includegraphics[width=0.35\textwidth]{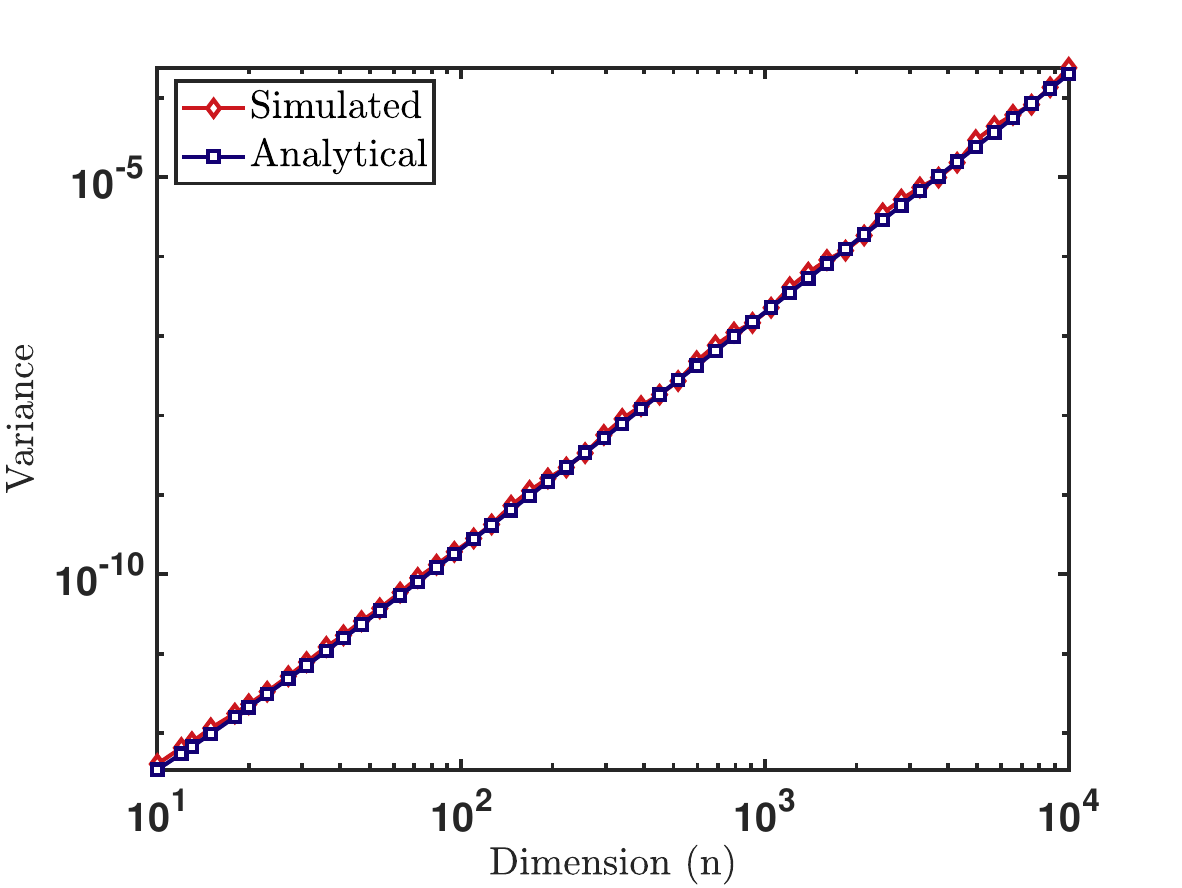}
    \label{fig:Inner_Product_Single_Nor_1_1}}
    \caption{Comparison between simulated variance and analytical variance, i.e., \eqref{eq:inner_e}, of the rounding error for the computation in single precision of the inner product with different input distribution.}
    \label{fig:Inner_Product}
\end{figure*}
In this section, we conduct a series of simulation experiments to validate the accuracy of our derivations presented in Section \ref{sec:app_la} and \ref{sec:wishart}.
\subsection{Simulation Setup}
The experiments are performed using MATLAB R2023b. While most computations are conducted in single precision, Section \ref{subsec:inner} employs $\mathtt{fp16}$ and $\mathtt{bfloat16}$ arithmetic. MATLAB $\mathtt{single.m}$ function simulates single precision, and the rounding function $\mathtt{chop.m}$ introduced in \cite{higham2019simulating} is used for simulating $\mathtt{fp16}$ and $\mathtt{bfloat16}$ arithmetic. Exact results for inner products and other matrix computations are obtained in double precision.

To ensure reproducibility, we initialize the random number generator with $\mathtt{rng}(1)$ at the start of each script generating a figure in this section. Each experiment is repeated $10000$ times for various problem dimensions $n$ and degrees of freedom $m$. Randomly generated matrices and vectors are employed, comparing different input distributions, including random uniform $\mathcal{U}\left(0,1\right)$, random uniform $\mathcal{U}\left(-1,1\right)$, random Gaussian $\mathcal{N}\left(0,1\right)$, random Gaussian $\mathcal{N}\left(1,1\right)$, and random chi-squared $\chi^2_m$ distribution.

Moreover, we use loops to implement the inner product and other matrix computations in MATLAB, where each operation involves rounding. This approach is necessary because MATLAB functions rely on the specifics of how the underlying BLAS operation is coded and optimized. For example, the accumulation of sums may involve extra precision for intermediate quantities, leading to inaccuracies in our results \cite{castaldo2009reducing}. Additionally, given the small value of $\sigma^2$ (approximately $10^{-16}$ for single precision), we employ the MATLAB $\mathtt{vpa.m}$ and $\mathtt{digital.m}$ functions to ensure accurate numerical computation of our analytical expressions, with $1000$ significant digits specified.

\subsection{General Rounding Error Analysis}
\subsubsection{Inner Products}
\label{subsec:inner}
We first show the numerical results for inner products $s= {\bf x}^T{\bf y}$ to validate the correctness of \eqref{eq:inner_var} and give some insights.

\textbf{Different input random vectors.}      
We consider the case where the entries of the input vectors have different distributions. Specifically, uniform $\mathcal{U}\left(0,1\right)$ and $\mathcal{U}\left(-1,1\right)$ distributions, as well as Gaussian $\mathcal{N}\left(0,1\right)$ and $\mathcal{N}\left(1,1\right)$ distributions. As shown in Fig. \ref{fig:Inner_Product}, the comparison between simulated variance and analytical variance, i.e., \eqref{eq:inner_e}, of the rounding error for the computation in single precision is illustrated. We can find that the analytical and simulated curves are very tight in different distributions, confirming our derived results' correctness. Moreover, the variance of the rounding error grows exponentially with the input dimension since the rounding errors accumulate along the vector dimension. Note that inputs with zero mean exhibit lower variances than those with nonzero mean. 

\textbf{Compared with other worst-case bounds.}
We employ the MSE, i.e., a classic statistical metric to analyze rounding error \cite{1056490}, of the computed results as our comparative metric against other worst-case bounds. Both the deterministic bounds \cite{higham2002accuracy,doi:10.1137/19M1270434} and the probabilistic bounds \cite{doi:10.1137/18M1226312,doi:10.1137/20M1314355,doi:10.1137/19M1270434} are considered. The metric can be mathematically formulated as:
\begin{equation}
    \mathbb{E}\left(\left| \hat{s} -s  \right|^2 \right) = \mathbb{E}\left(\left|\Delta s  \right|^2 \right).
\end{equation}
From {Theorem} \ref{the:inner}, we have 
\begin{align}
    \mathbb{E}\left(\left| \hat{s} -s  \right|^2 \right) &=\mathbb{E}\left(\Delta s ^2 \right)=\mathbb{V}\left(\Delta s   \right)\nonumber\\
    &\approx\hbar\left(\mu_x,\sigma_x,\mu_y,\sigma_y,n,\sigma \right).\label{eq:mse_r}
\end{align}
\begin{figure}[t]
    \centering
    \subfloat[Random uniform $\mathcal{U}\left(0,1\right)$ vectors.]{\includegraphics[width=0.35\textwidth]{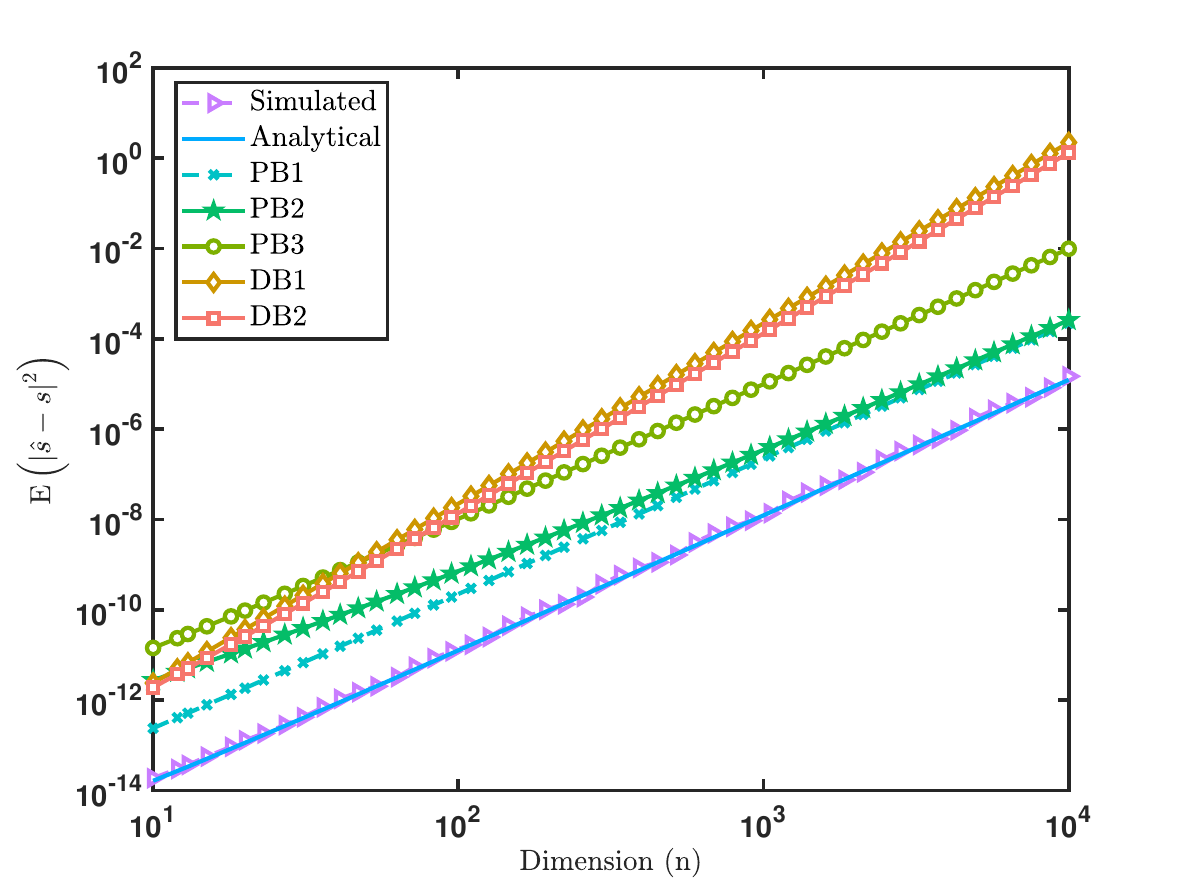}
    \label{fig:Compare}}
    \vfill
    \subfloat[Random uniform $\mathcal{U}\left(-1,1\right)$ vectors.]{\includegraphics[width=0.35\textwidth]{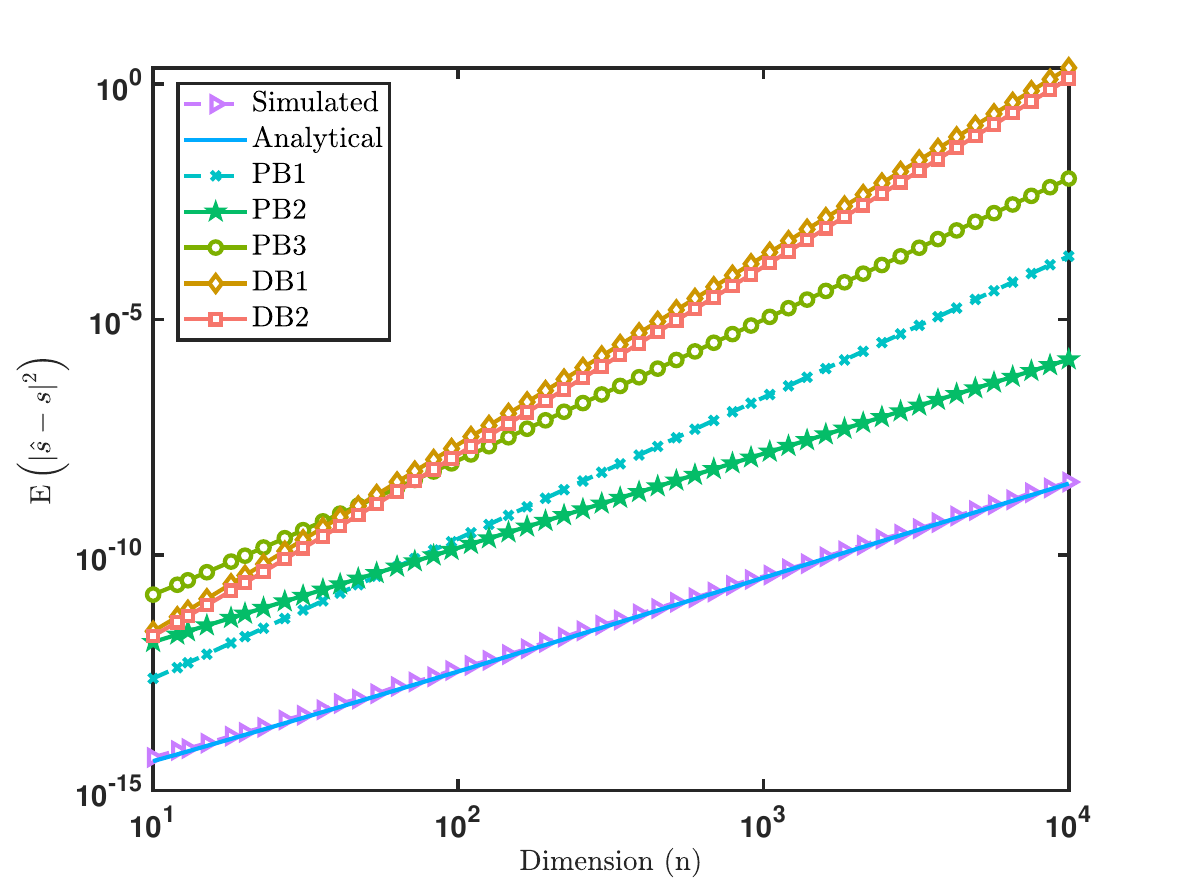}
    \label{fig:Compare_Uni}}
    \caption{Comparison between the analytical results and other worst-case bounds for the computation in single precision of the inner product with different input distributions. Here, $\lambda = 1$ and $\zeta = 10^{-16}$.}
    \label{fig:Compare_worst_case}
\end{figure}
Then, for the worst-case bounds and using \cite[Theorem 3.1, Eq. (3.6)]{doi:10.1137/18M1226312}, we have
\begin{align}
    ({\rm DB1})~\mathbb{E}\left(\left| \hat{s} -s  \right|^2 \right) 
     &\leq \gamma_n^2  \mathbb{E}\left(\left(\left| {\bf x}\right|^T \left|{\bf y}  \right|\right)^2 \right) \label{eq:mse_worst},\\
     ({\rm PB1})~\mathbb{E}\left(\left| \hat{s} -s  \right|^2 \right) 
     &\leq \gamma_n^2(\lambda)  \mathbb{E}\left(\left(\left| {\bf x}\right|^T \left|{\bf y}  \right|\right)^2  \right) \label{eq:mse_worst_p},
\end{align}
where $\gamma_n = \frac{nu}{1-nu}$ and $\gamma_n^2(\lambda) = \exp \left( \lambda \sqrt{n}u+\frac{nu^2}{1-nu} \right) -1$. $\lambda$ is a positive constant that can be freely chosen and controls the probability of failure of the bound, which is a monotonically decreasing function of $\lambda$.

Next, based on the probabilistic bound in \cite[Theorem 3.2]{doi:10.1137/20M1314355}, we obtain

{\footnotesize\begin{align}
   \mathbb{E}\left(\left| \hat{s} -s  \right|^2 \right) &\leq \mathbb{E}\left(\left| \left(\lambda\left|\mu_x\mu_y\right|n^{\frac{3}{2}} + \left(\lambda^2 + 1\right)C_xC_y n \right)u+\mathcal{O}\left(u^2\right) \right|^2 \right) \nonumber\\
   ({\rm PB2})~ &\approx \left(\lambda\left|\mu_x\mu_y\right|n^{\frac{3}{2}} + \left(\lambda^2 + 1\right)C_xC_y n \right)^2u^2 \label{eq:mse_shaper},
\end{align}}
where $\left|x_i\right|\leq C_x$ and $\left|y_i\right|\leq C_y$ for $i=1:n$. For example, given $x_i \sim \mathcal{U}(0,1)$, we have $C_x = 1$.

Finally, we utilize \cite[Theorem 3.2 \& 3.3]{doi:10.1137/19M1270434}, yielding
\begin{align}
    ({\rm DB2})~\mathbb{E}\left(\left| \hat{s} -s  \right|^2 \right) & \leq n\mathbb{E}\left( \sum_{k = 1}^n c_k^2\right)\label{eq:det_2},\\
    ({\rm PB3})~\mathbb{E}\left(\left| \hat{s} -s  \right|^2 \right) & \leq 2\ln\left(2/\zeta \right)\mathbb{E}\left( \sum_{k = 1}^n c_k^2\right)\label{eq:pro_2},
\end{align}
where $c_1 = \left|x_1y_1\right|\beta_n$, $c_k = \left|x_ky_k\right|\beta_{n-k+2}$, $\beta_k = \left(1+u\right)^k - 1,2\leq k\leq n$ and $\zeta$ is the failure probability of \eqref{eq:pro_2}.

\begin{figure}[t]
    \centering
    \subfloat[$\mathtt{fp16}$ arithmetic.]{\includegraphics[width=0.35\textwidth]{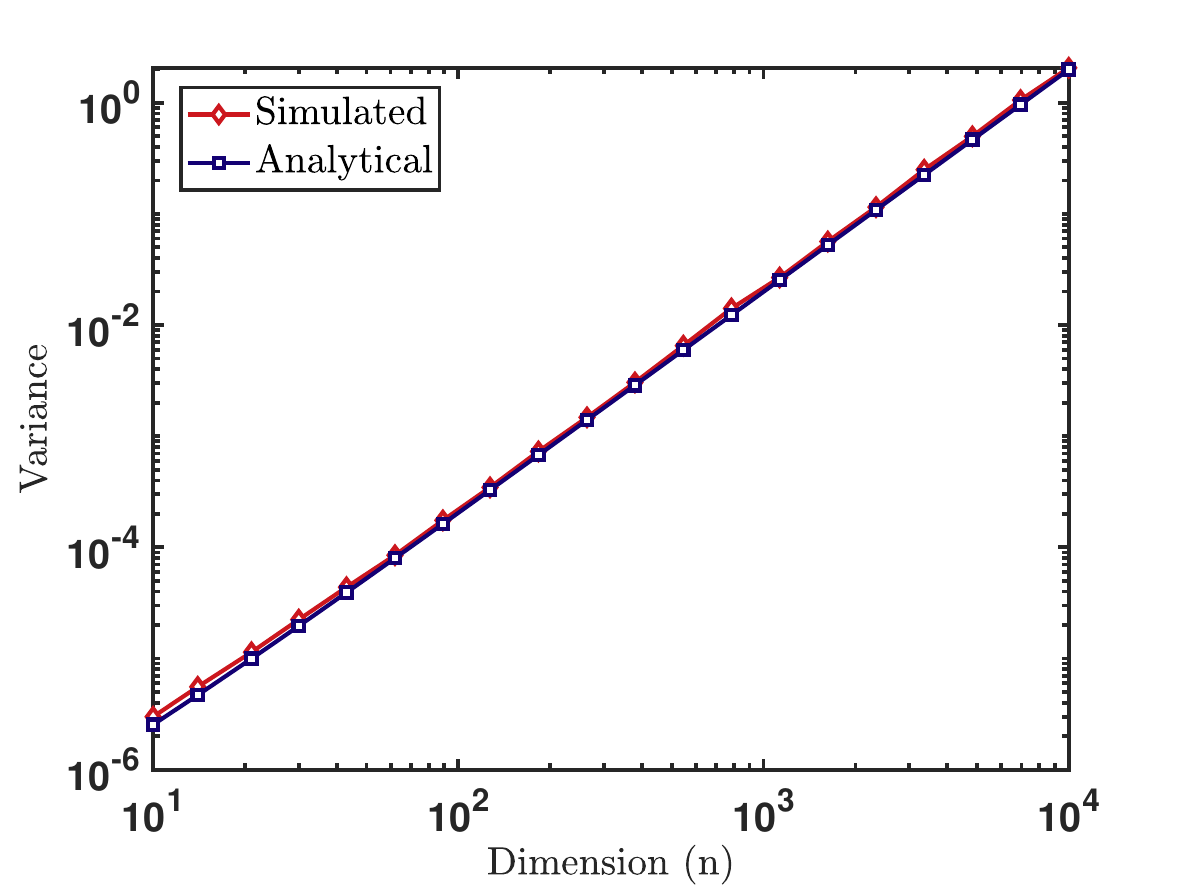}
    \label{fig:fp16}}
    \vfill
    \subfloat[$\mathtt{bfloat16}$ arithmetic.]{\includegraphics[width=0.35\textwidth]{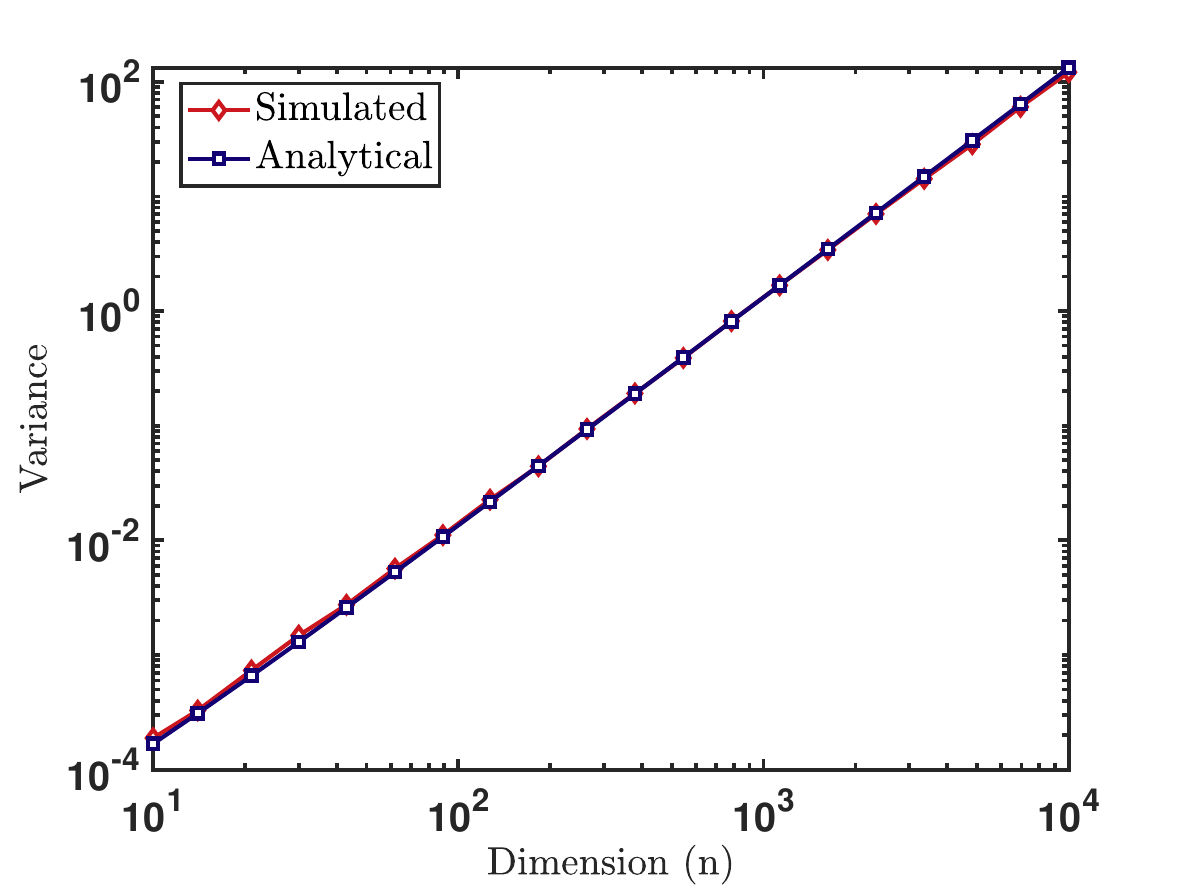}
    \label{fig:bfp16}}
    \caption{Comparison between simulated variance and analytical variance of the rounding error for the computation in lower precision of the inner product with different input distribution with random Gaussian $\mathcal{N}\left(0,1\right)$ vectors.}
    \label{fig:Lower}
\end{figure}

Therefore, using \eqref{eq:mse_r}, \eqref{eq:mse_worst}, \eqref{eq:mse_worst_p}, \eqref{eq:mse_shaper}, \eqref{eq:det_2} and \eqref{eq:pro_2}, we compare the analytical outcomes with other worst-case bounds for computing the inner product with random uniform $\mathcal{U}\left(0,1\right)$ and $\mathcal{U}\left(-1,1\right)$ vectors in single precision, as illustrated in Fig. \ref{fig:Compare} and \ref{fig:Compare_Uni}, respectively. It's clear that our analytical results align much more closely with the actual MSE than two deterministic and three probabilistic (We set $\lambda = 1$  and $\zeta = 10^{-16}$, which are classic parameters suggested in \cite{doi:10.1137/18M1226312,doi:10.1137/19M1270434}.) bounds for MSE. Specifically, the analytical expression tends to be at least two orders of magnitude tighter than other worst-case bounds. Notably, most of the bounds are even not approximate close-form expressions.


\textbf{Lower precision arithmetic.}
Now we repeat the experiment conducted in Section \ref{subsec:inner} using Gaussian $\mathcal{N}\left(0,1\right)$ distribution, but this time employing precision lower than single precision to compute the inner product $s= {\bf x}^T{\bf y}$. The simulation results using $\mathtt{fp16}$ and $\mathtt{bfloat16}$ arithmetic are depicted in Fig. \ref{fig:fp16} and Fig. \ref{fig:bfp16}, respectively. Notably, our analytical results remain highly consistent with the simulated curves even at lower precision, affirming the accuracy of our derivation.

\subsubsection{A Case Where {Model} \ref{model:pro} is Invalid}
\begin{figure}[t]
    \centering
    \subfloat[Simulated and analytical variance at each loop $i$ with random vectors with $x_i\sim \mathcal{N}\left(0,1\right)$ and $y_i = x_ih,h\sim\mathcal{N}\left(0,1\right)$.]{\includegraphics[width=0.35\textwidth]{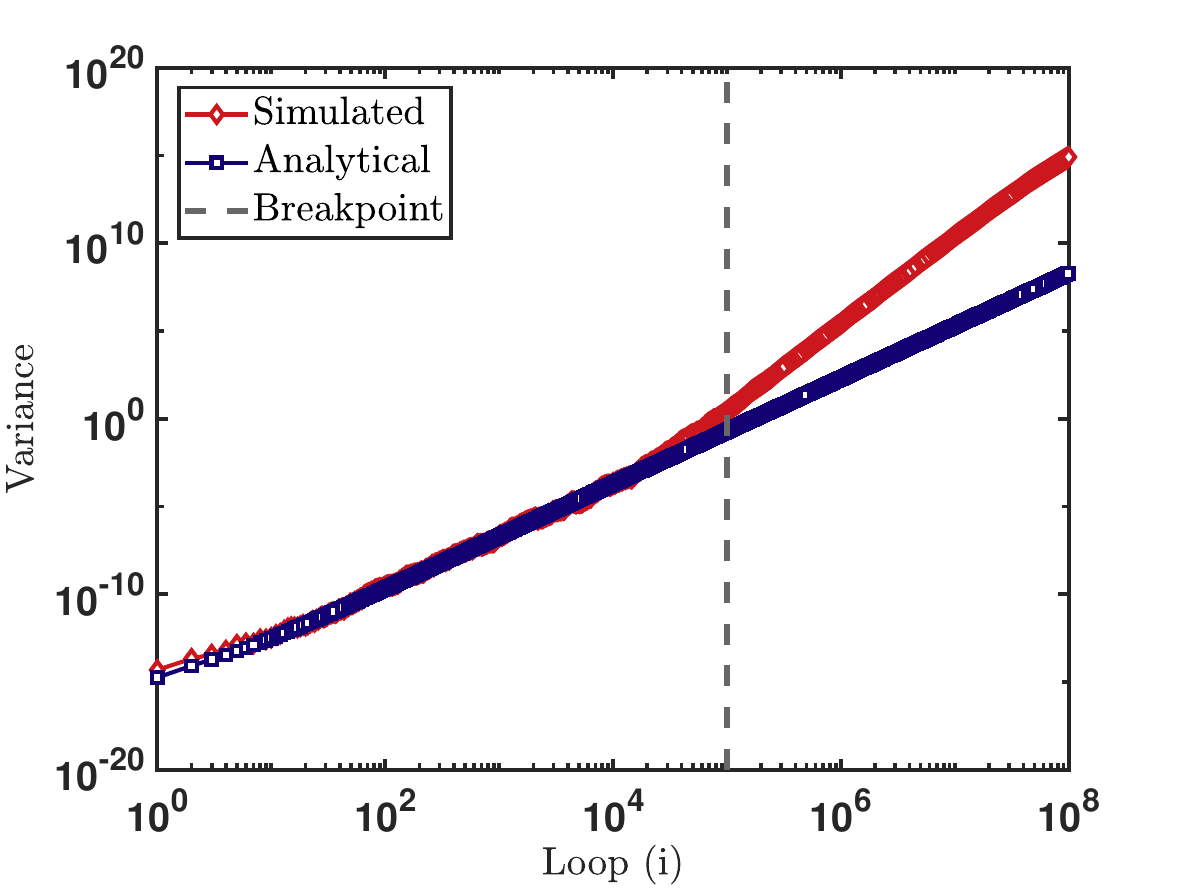}
    \label{fig:Inner_Large}}
    \vfill
    \subfloat[Distribution of $\delta_i$ for $i=1:10^5$ (top) and $i = 10^5:10^8$ (bottom), where the red curve is from $\eqref{eq:pdf_delta}$ and the purple curve is simulated.]{\includegraphics[width=0.35\textwidth]{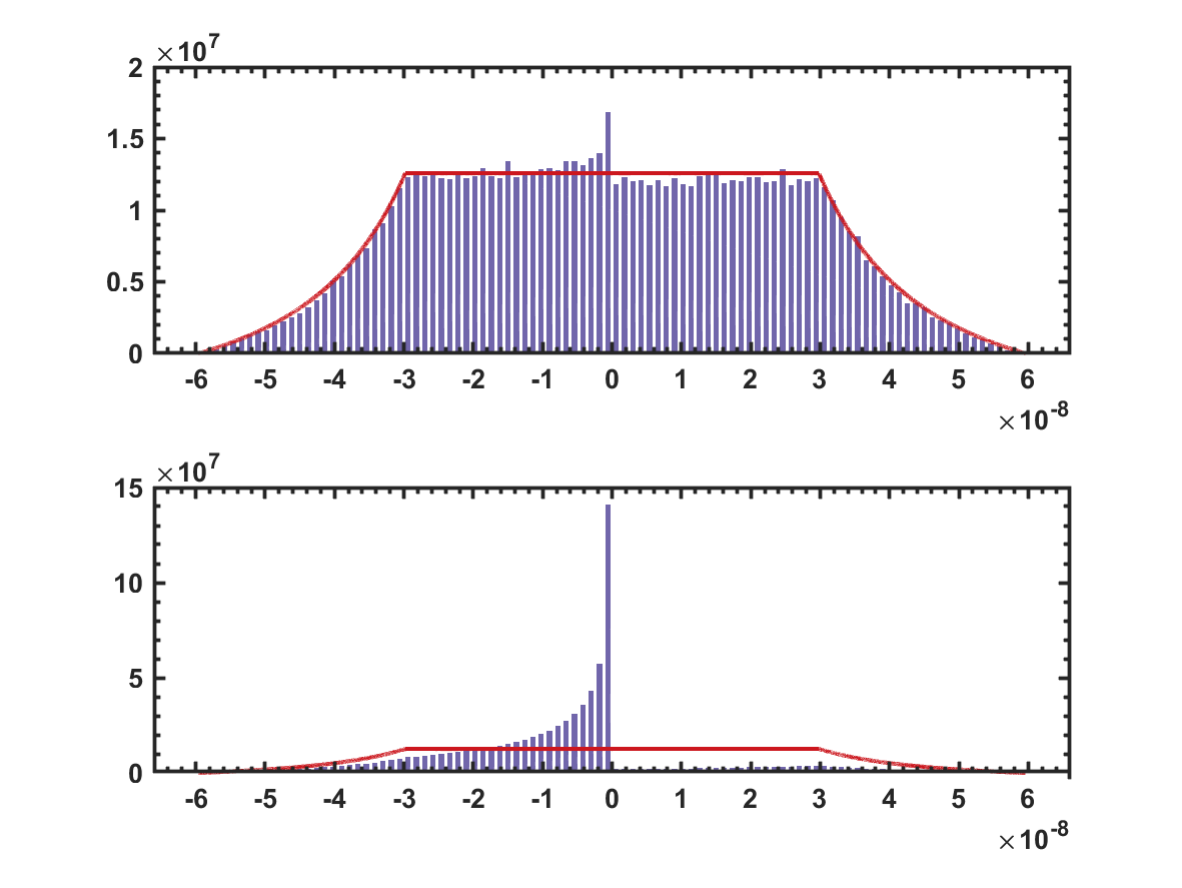}
    \label{fig:distribution}}
    \caption{Computation in single precision of the inner product $s={\bf x}^T{\bf y}$ with dependent random vectors of dimension $n = 10^8$.}
    \label{fig:model_invaild}
\end{figure}

In this subsection, we explore a scenario where {Model} \ref{model:pro} fails and the analytical results are inaccurate, using the computation of the inner product $s= {\bf x}^T{\bf y}$. Specifically, we consider {large dependent} input vectors of dimension $n = 10^8$ with $x_i\sim \mathcal{N}\left(0,1\right)$ and $y_i = x_ih,h\sim\mathcal{N}\left(0,1\right)$. As illustrated in Fig. \ref{fig:Inner_Large}, we depict the simulated and analytical variances at each loop iteration $i$ of the inner products. Notably, while the analytical and simulated curves closely match for $i=1:10^5$, a significant gap emerges for $i>10^5$, indicating the invalidity of {Model} \ref{model:pro} in this context. This discrepancy arises because the distribution of the relative error $\delta$ changes with the large dependent input vectors. Although \eqref{eq:pdf_delta} in {Model} \ref{model:pro} provides a reasonable approximation for $i=1:10^5$, it does not hold for $i>10^5$, as evidenced by Fig. \ref{fig:distribution}. 

\subsubsection{Matrix-vector and matrix-matrix products}
\begin{figure}[t]
    \centering
    \subfloat[${\bf R}_{\Delta {\bf C}}(2,2)$ as a function of $n$. Here, $m = 10$ and $p = 10$.]{\includegraphics[width=0.3\textwidth]{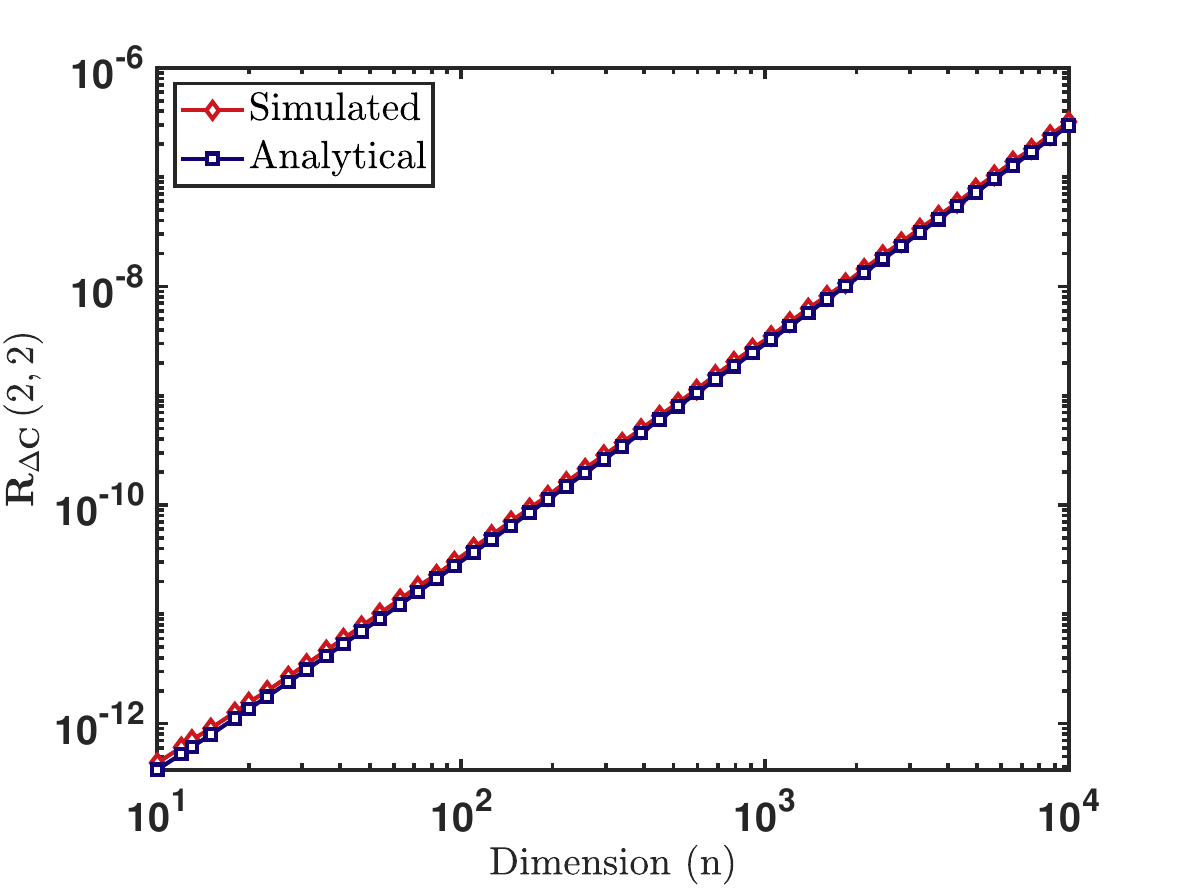}
    \label{fig:MM_n}}
    \vfill
    \subfloat[${\bf R}_{\Delta {\bf C}}(2,2)$ as a function of $p$. Here, $m = 10$ and $n = 10$.]{\includegraphics[width=0.3\textwidth]{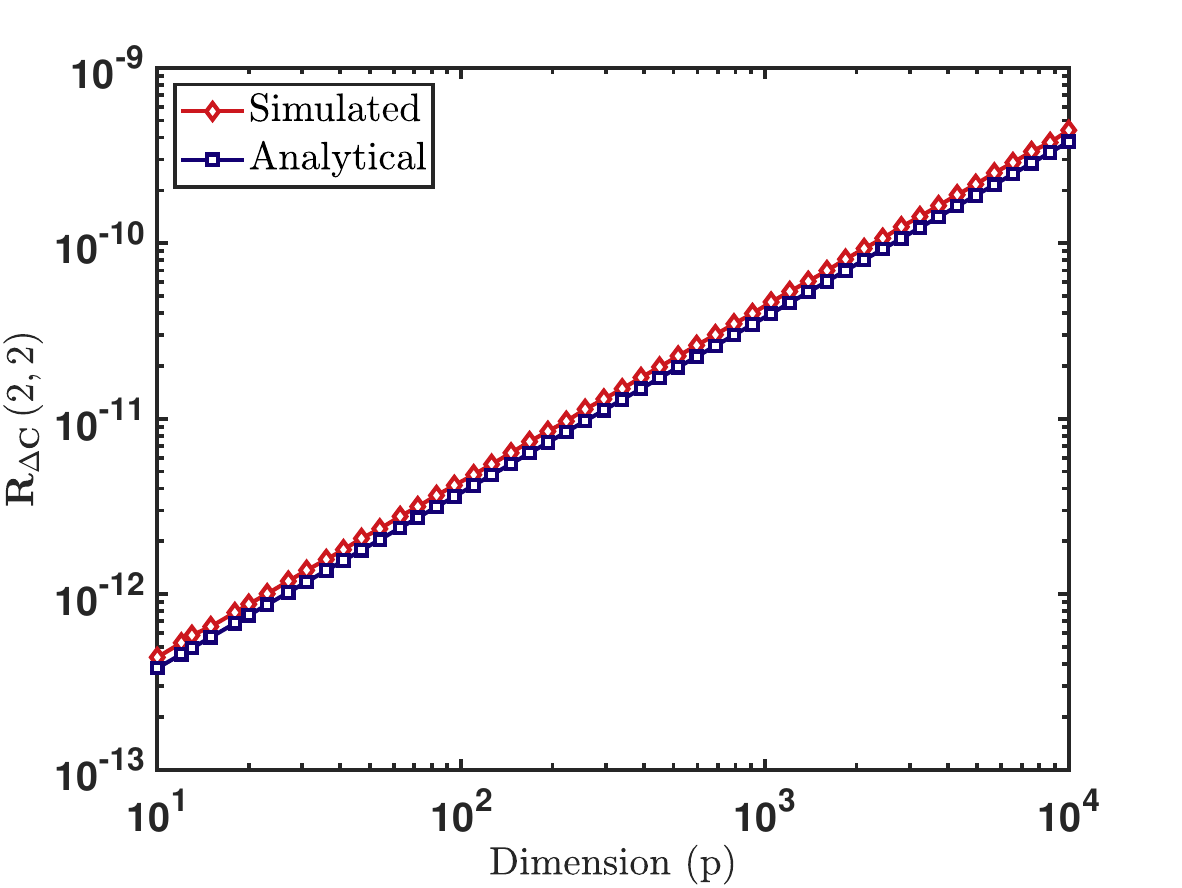}
    \label{fig:MM_p}}
    \vfill
    \subfloat[${\bf R}_{\Delta {\bf C}}(2,2)$ as a function of $m$. Here, $n = 10$ and $p = 10$.]{\includegraphics[width=0.3\textwidth]{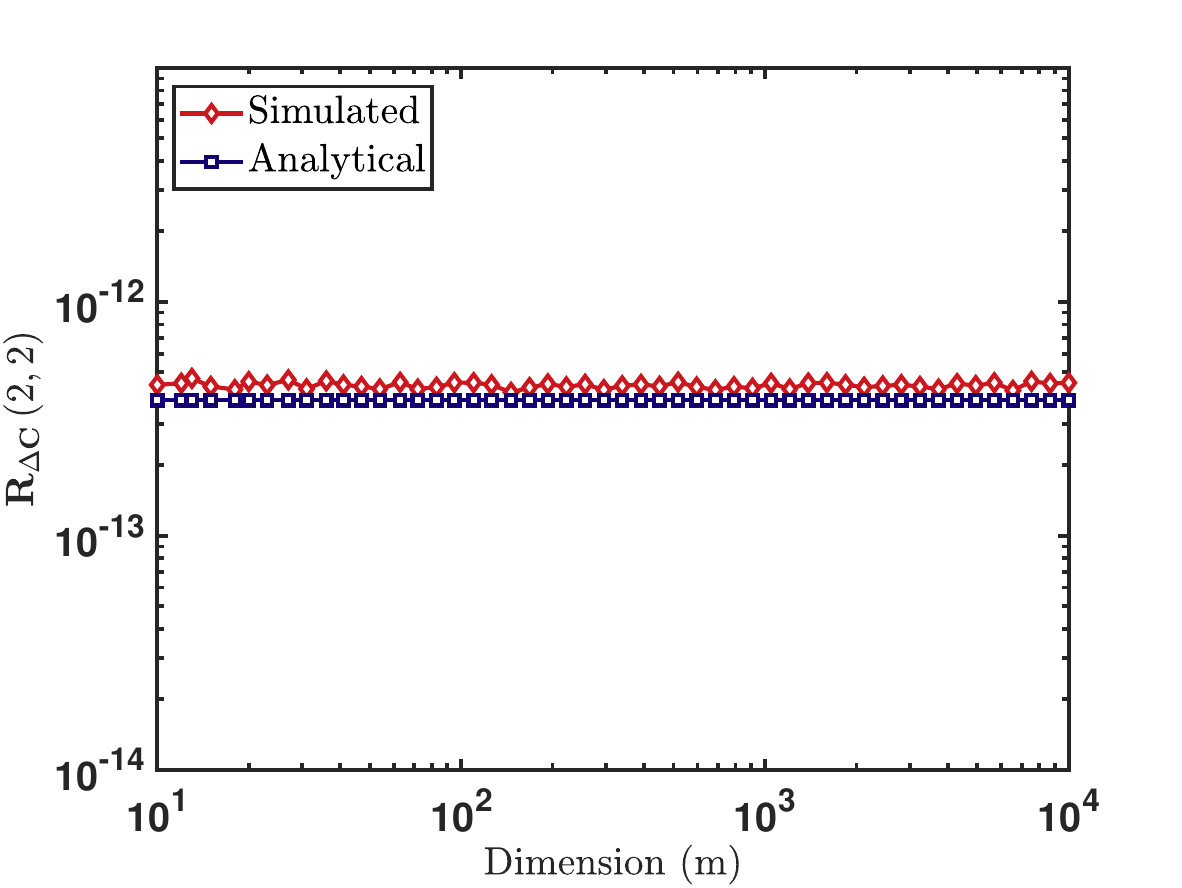}
    \label{fig:MM_m}}
    \caption{Comparison between simulated autocorrelation matrix and analytical autocorrelation matrix, i.e., \eqref{eq:mm_var}, of the rounding error for the computation in single precision of the matrix-matrix product with different dimensions using the second-row second-column element ${\bf R}_{\Delta {\bf C}}(2,2)$ as an example.}
    \label{fig:MM_Product}
\end{figure}

This subsection considers the computation of matrix-vector and matrix-matrix products in single precision. Using the matrix-matrix product $\bf C=AB$ as a case study, where ${\bf A}\in \mathbb{R}^{m\times n}$ and ${\bf B}\in \mathbb{R}^{n\times p}$, we plot the simulated and analytical curves for the autocorrelation matrix element ${\bf R}_{\Delta {\bf C}}(2,2)$, demonstrating the impact of dimension $m$, $n$, and $p$ as depicted in Fig. \ref{fig:MM_Product}. We can observe that the analytical results accurately anticipate variations in different dimensions, affirming the correctness of our derived outcomes. Additionally, the autocorrelation matrix of the rounding error remains unaffected by $m$, as illustrated in Fig. \ref{fig:MM_m}.

\subsection{Specific Rounding Error Analysis for Wishart Matrices}
\subsubsection{Triangular Systems}
\label{subsec:tri}
\begin{figure}[t]
    \centering
    \subfloat[Variance of $\Delta x_3$ as a function of DoF $m$. Here, $n = 5$.]{\includegraphics[width=0.35\textwidth]{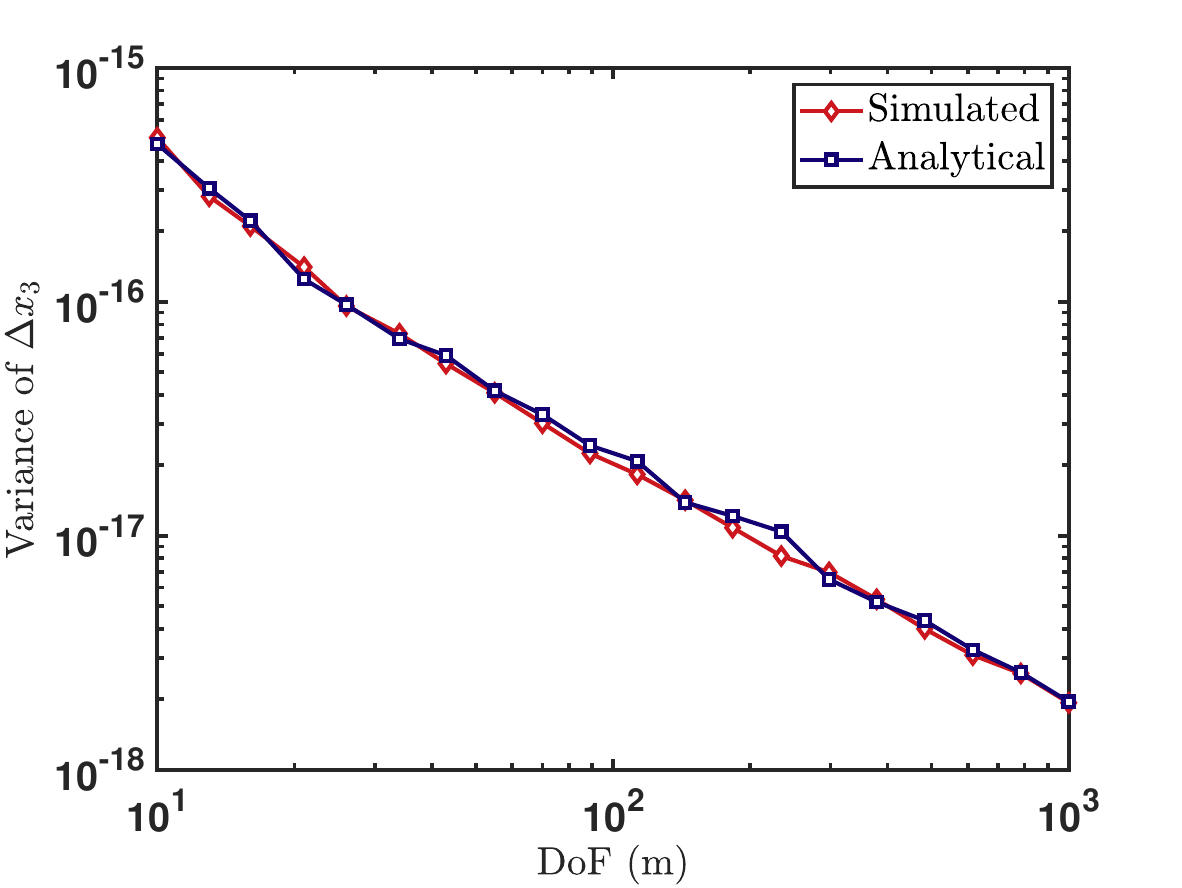}
    \label{fig:Tri_DoF}}
    \vfill
    \subfloat[Variance of $\Delta x_3$ as a function of dimension $n$. Here, $m = 1050$.]{\includegraphics[width=0.35\textwidth]{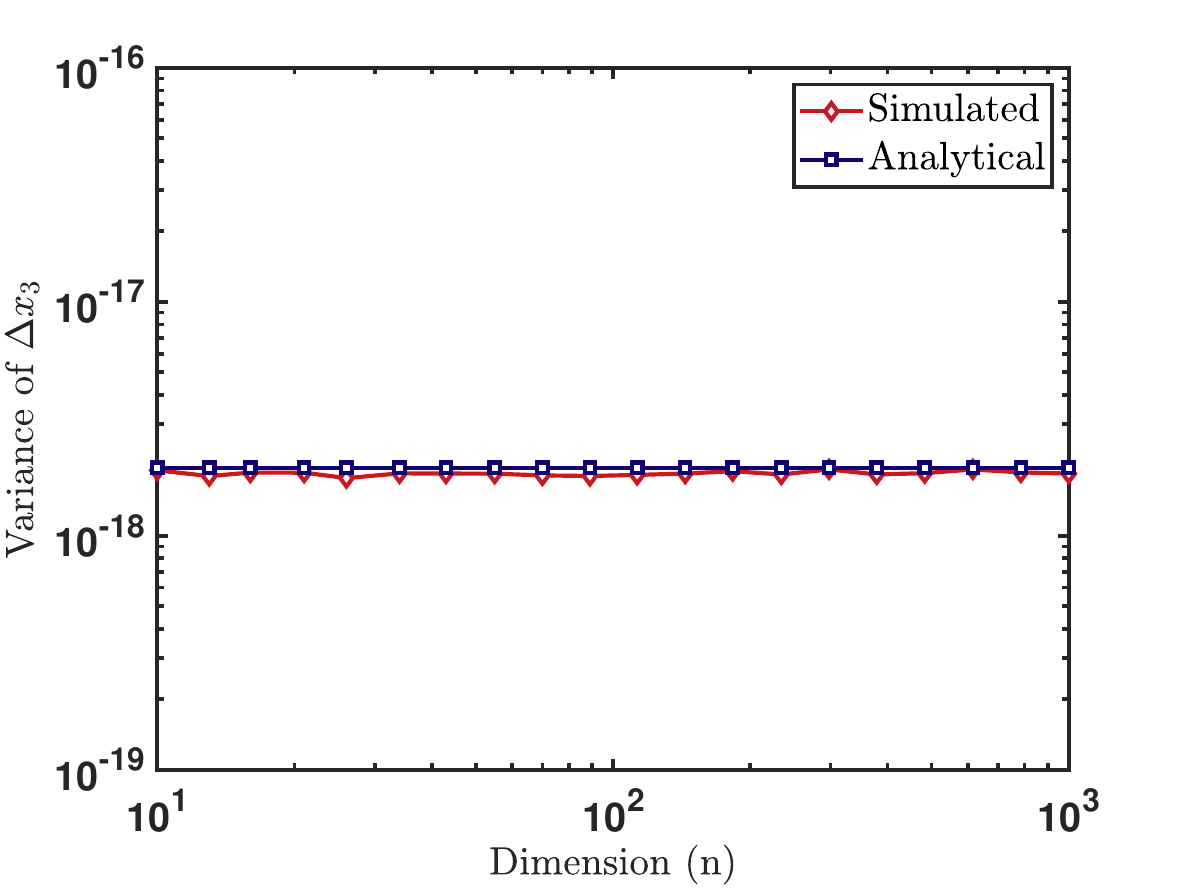}
    \label{fig:Tri_Dimen}}
    \caption{Comparison between simulated variance and analytical variance of the rounding error for the solution of triangular systems in single precision using $x_3$ as an example.}
    \label{fig:tri_simu}
\end{figure}

Next, in this subsection, we present the solution of triangular systems ${\bf Tx=b}$. For simplification, the variance of the rounding error for $x_3$ is shown as an example. From {Theorem} \ref{the:tri}, we have

{\footnotesize\begin{align}
    \mathbb{V}\left( \Delta x_3 \right) &\approx\frac{\left( 1+\sigma ^2 \right) ^3+\sum_{j=1}^{2}{\mathbb{V}\left( x_j \right) \left( 1+\sigma _{\psi _j}^{2} \right) \left( 1+\sigma ^2 \right) ^{i-j+2}}}{m-4}\nonumber\\
    &-\mathbb{V}\left( x_3 \right)\nonumber\\
    & = \frac{\left( 1+\sigma ^2 \right) ^3+\frac{\left( 1+\sigma ^2 \right) ^5}{m-2}+\frac{\left( m-1 \right) \left( 1+\sigma _{\psi _2}^{2} \right) \left( 1+\sigma ^2 \right) ^3}{\left( m-2 \right) \left( m-3 \right)}}{m-4}\nonumber\\
    &-\mathbb{V}\left( x_3 \right) \label{eq:tri_3},
\end{align}}
where 
\begin{align*}
    \sigma _{\psi _2}^{2} &= \frac{\left[ \left( 1+\sigma ^2 \right) ^2-1 \right] \left[ m+\left( 1+\sigma ^2 \right) ^2-1 \right]}{m-1},\\
    \mathbb{V}\left( x_3 \right) & = \frac{\left( m-2 \right) \left( m-3 \right) +2m-4}{\left( m-2 \right) \left( m-3 \right) \left( m-4 \right)}.
\end{align*}

In Fig. \ref{fig:tri_simu}, we compare simulated variance and analytical variance of the rounding error for the solution of triangular systems in single precision by \eqref{eq:tri_3}. Similar to the case for inner products and matrix-matrix products, the analytical and simulated curve is still very tight. Moreover, it is observed that the variance of the rounding error is independent of the input dimension and mainly depends on DoF $m$. 

\subsubsection{LU Factorization}
\begin{figure*}[t]
    \centering
    \subfloat[Variance of $\Delta u_{33}$ as a function of DoF $m$. Here, $n = 5$.]{\includegraphics[width=0.3\textwidth]{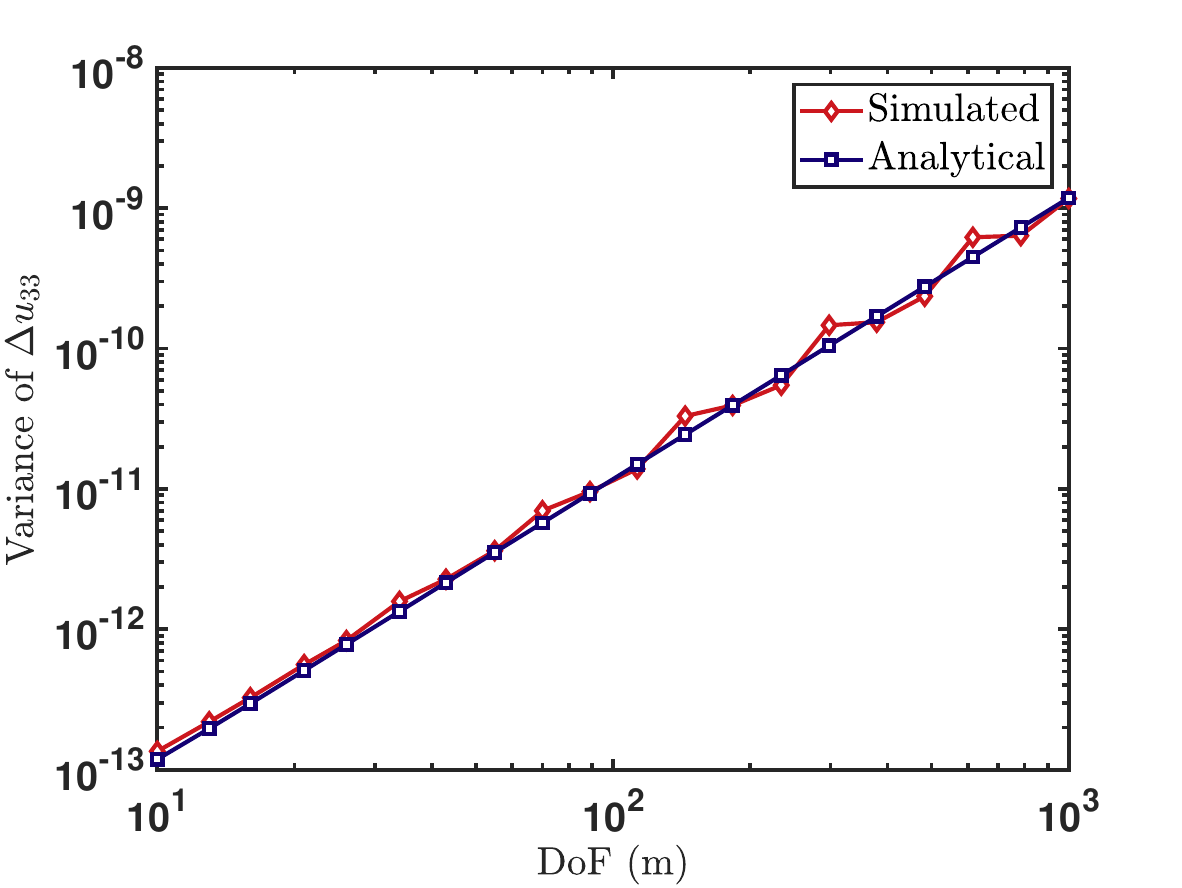}
    \label{fig:u_33_m}}
    \hfill
    \subfloat[Variance of $\Delta u_{35}$ as a function of DoF $m$. Here, $n = 5$.]{\includegraphics[width=0.3\textwidth]{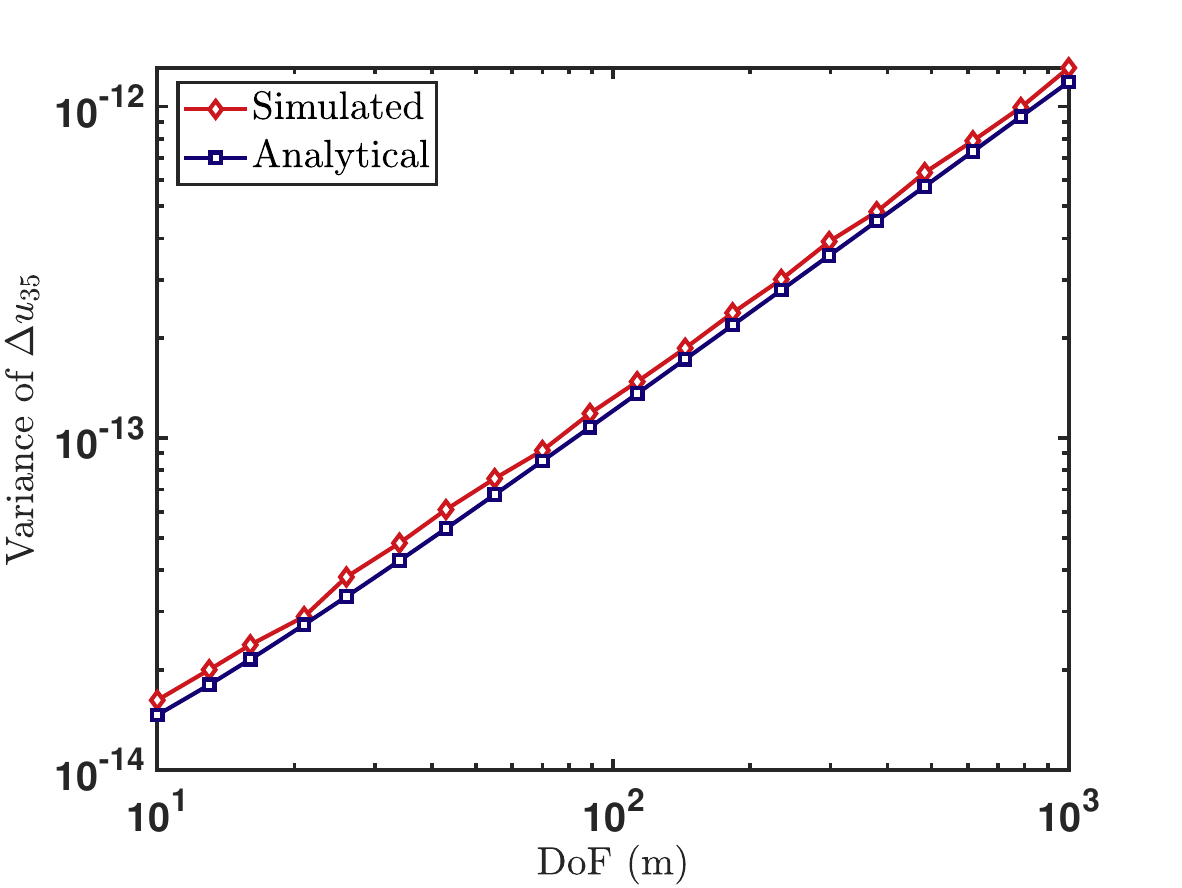}
    \label{fig:u_35_m}}
    \hfill
    \subfloat[Variance of $\Delta l_{43}$ as a function of DoF $m$. Here, $n = 5$.]{\includegraphics[width=0.3\textwidth]{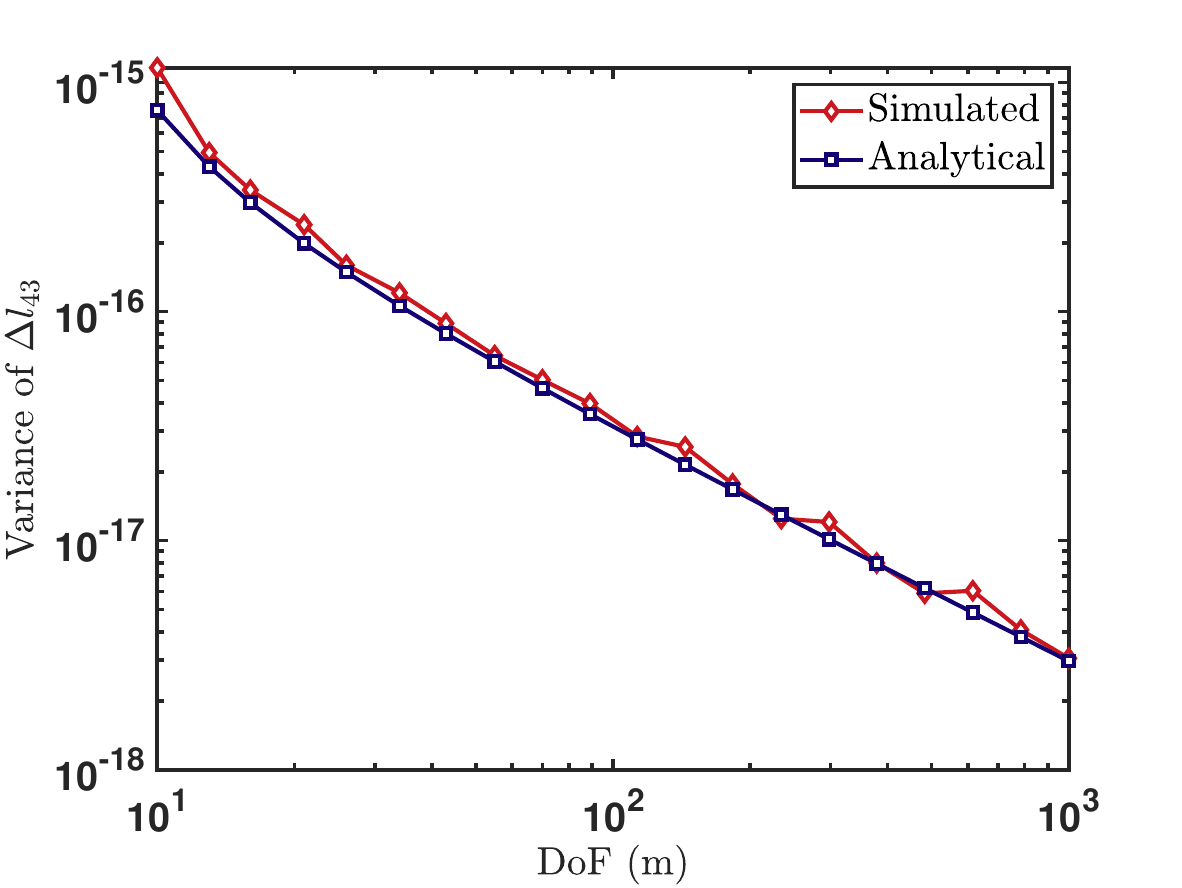}
    \label{fig:l_43_m}}\\
    \subfloat[Variance of $\Delta u_{33}$ as a function of dimension $n$. Here, $m = 1050$.]{\includegraphics[width=0.3\textwidth]{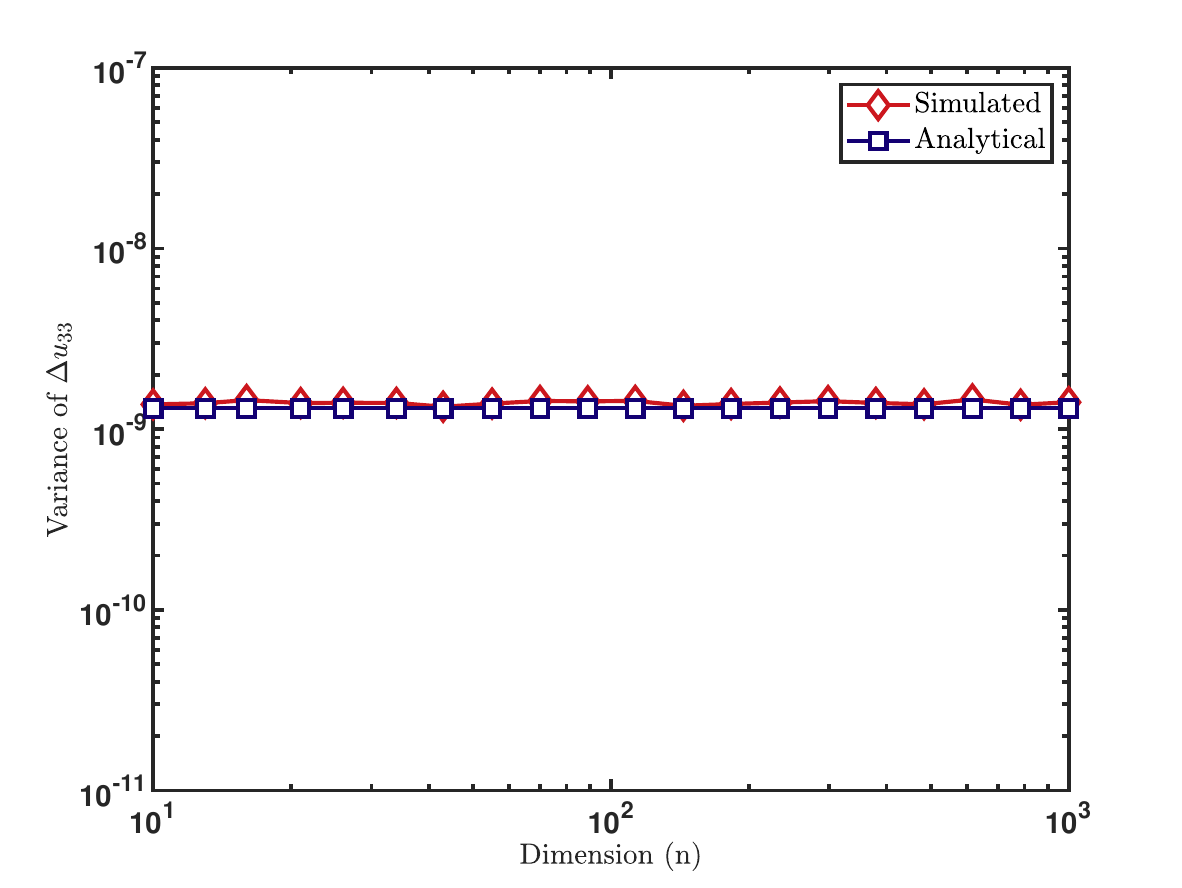}
    \label{fig:u_33_n}}
    \hfill
    \subfloat[Variance of $\Delta u_{35}$ as a function of dimension $n$. Here, $m = 1050$.]{\includegraphics[width=0.3\textwidth]{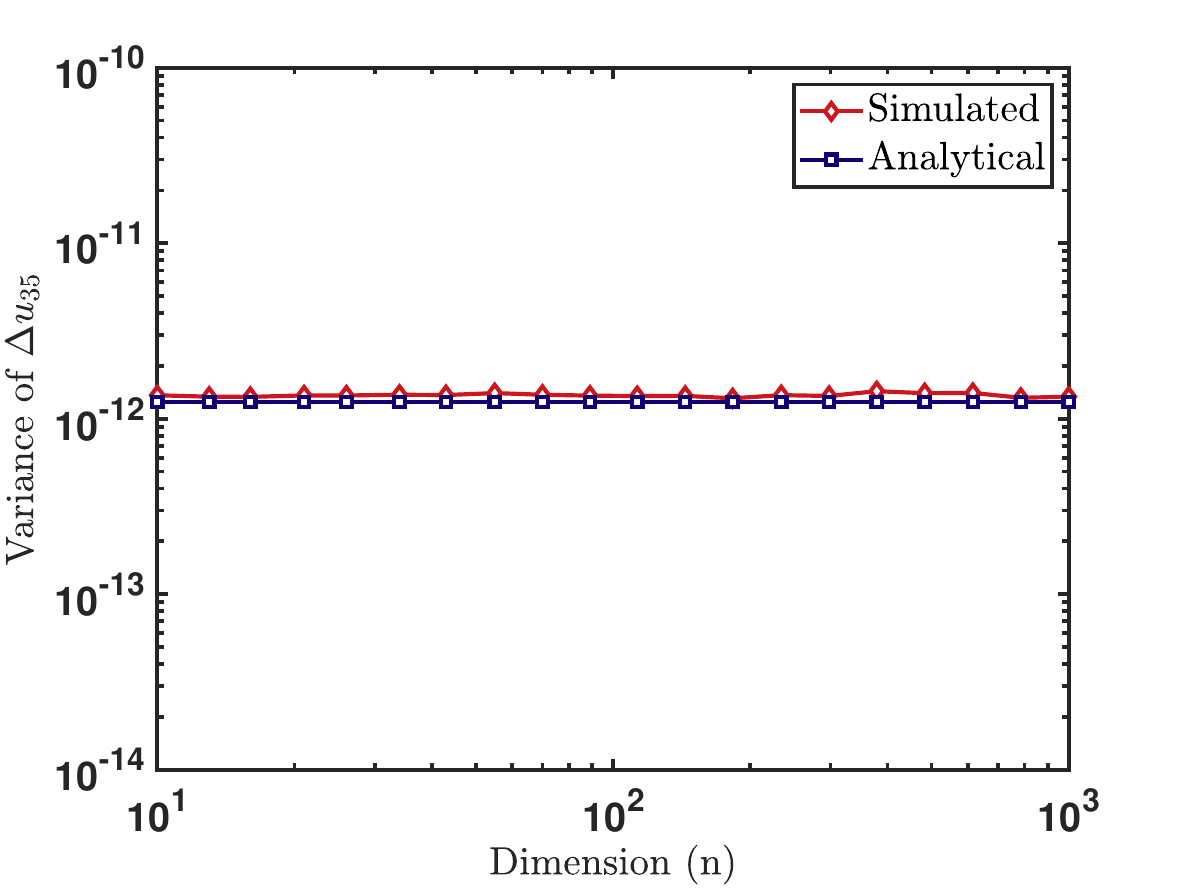}
    \label{fig:u_35_n}}
    \hfill
    \subfloat[Variance of $\Delta l_{43}$ as a function of dimension $n$. Here, $m = 1050$.]{\includegraphics[width=0.3\textwidth]{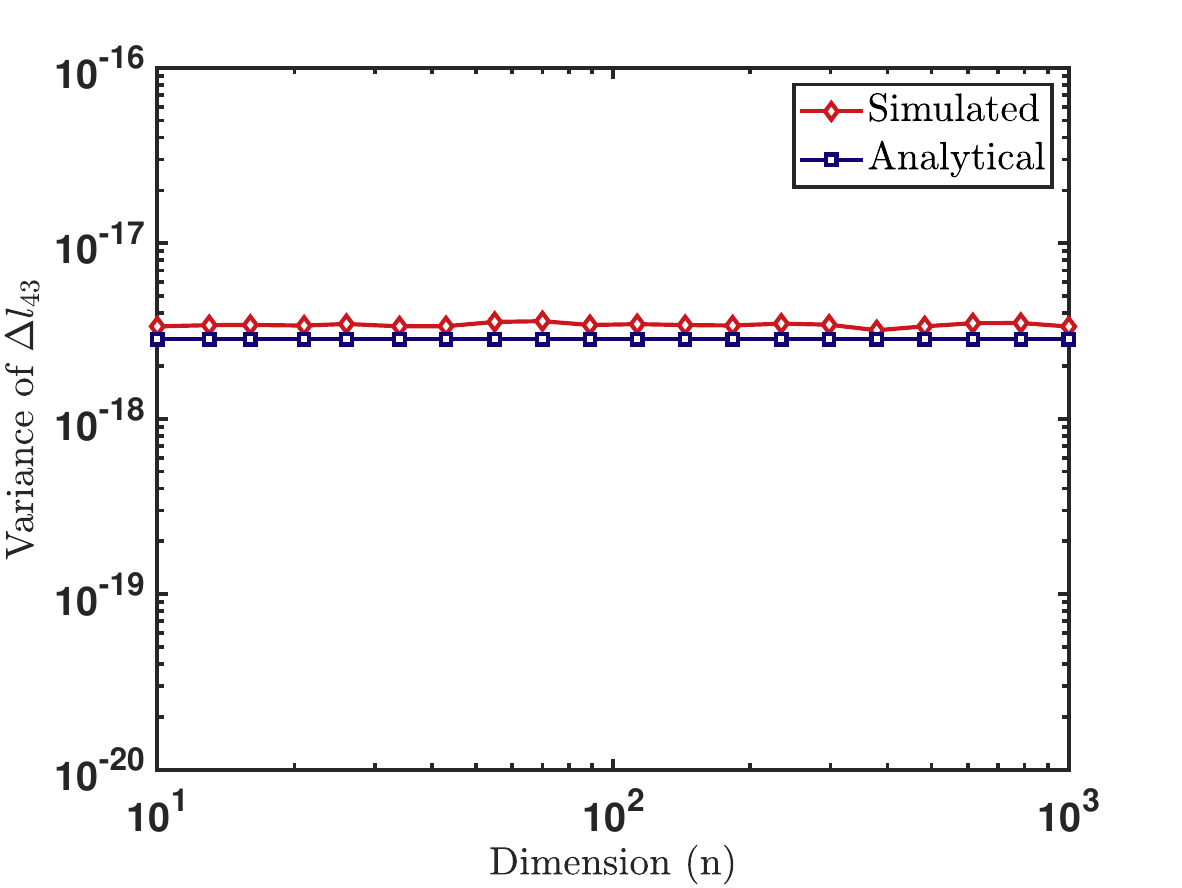}
    \label{fig:l_43_n}}
    \caption{Comparison between simulated variance and analytical variance of the rounding error for the computation in single precision of LU factorization using $u_{33}$, $u_{53}$ and $l_{43}$ as an example.}
    \label{fig:LU}
\end{figure*}

Finally, we consider the rounding error of LU factorization $\bf A =LU$. For simplification, the variances of the rounding errors for $u_{33}$, $u_{35}$ and $l_{43}$ are shown as an example. Using {Theorem} \ref{the:lu} with $k = 3$, we can obtain

{\footnotesize\begin{align}
    \mathbb{V}\left( \Delta u_{33} \right) &\approx \left( m^2-4 \right) \left[ \left( 1+\sigma ^2 \right) ^2-1 \right]\nonumber \\
    &+3\left( 1+\sigma ^2 \right) ^2\left[ \left( 1+\sigma _{\epsilon _2}^{2} \right) \left( 1+\sigma _{\eta _2}^{2} \right) +\left( 1+\sigma ^2 \right) ^2 \right]  \nonumber\\
    &-2\left( m+2 \right) \sigma ^2-6, \label{eq:u33}\\
    \mathbb{V}\left( \Delta u_{3j} \right) &\approx \left( m-2 \right) \left[ \left( 1+\sigma ^2 \right) ^2-1 \right] +\left( 1+\sigma ^2 \right) ^4\nonumber\\
    &+\left( 1+\sigma _{\epsilon _2}^{2} \right) \left( 1+\sigma _{\eta _2}^{2} \right) \left( 1+\sigma ^2 \right) ^2 \nonumber \\
    &-2\left( 1+\sigma ^2 \right), \,\,\,\,\,\,  j \neq 3,\label{eq:u35} \\
    \mathbb{V} \left( \Delta l_{i3} \right) &\approx \frac{\left( 1+\sigma _{\eta _3}^{2} \right) \left( 1+\sigma ^2 \right) ^3-1}{m-6}
    \nonumber\\
    &+\frac{\left( 1+\sigma _{\eta _3}^{2} \right) \left[ \left( 1+\sigma ^2 \right) ^5+\left( 1+\sigma _{\epsilon _2}^{2} \right) \left( 1+\sigma _{\eta _2}^{2} \right) \left( 1+\sigma ^2 \right) ^3 \right]}{\left( m-4 \right) \left( m-6 \right)}
    \nonumber\\
    &-\frac{2\left( 1+\sigma _{\eta _3}^{2} \right) \left( 1+\sigma ^2 \right) ^2}{\left( m-4 \right) \left( m-6 \right)}\label{eq:l43},
\end{align}}
where 
{\small \begin{align*}
    \sigma _{\epsilon _2}^{2}&=\frac{\left( m-6 \right) \left[ 1+\frac{\left( m^2-4 \right) \sigma ^2+3\left[ \left( 1+\sigma ^2 \right) ^3-1 \right]}{m^2-1} \right]}{m-5}\nonumber\\
    &+\frac{\frac{\left( m^2-4 \right) \sigma ^2+3\left[ \left( 1+\sigma ^2 \right) ^3-1 \right]}{m^2-1}\left( 1+\sigma ^2 \right) ^4-1}{m-5},\nonumber\\
    \sigma _{\eta _2}^{2}&=\frac{\left( m-2 \right) \sigma ^2+\left( 1+\sigma ^2 \right) ^3-1}{m-1},\nonumber\\
    \sigma _{\eta _3}^{2}&=\frac{\left( m^2-4 \right) \left[ \left( 1+\sigma ^2 \right) ^2-1 \right]}{m^2-2m}-\frac{2\left( m+2 \right) \sigma ^2}{m^2-2m}\nonumber\\
    &+\frac{3\left\{ \left( 1+\sigma ^2 \right) ^2\left[ \left( 1+\sigma _{\epsilon _2}^{2} \right) \left( 1+\sigma _{\eta _2}^{2} \right) +\left( 1+\sigma ^2 \right) ^2 \right] -2 \right\}}{m^2-2m}.
\end{align*}}

Then, we compare the simulated and analytical variance of the rounding error for solving triangular systems in single precision using equations \eqref{eq:u33}, \eqref{eq:u35}, and \eqref{eq:l43} in Fig. \ref{fig:LU}. Consistent with the findings in \ref{subsec:tri}, our analytical results accurately predict the trend of the actual variance, validating our derivation. Furthermore, the variance of the rounding error depends solely on the DoF and precision, remaining uncorrelated with the dimension of ${\bf A}$.


\section{Conclusions}
\label{sec:conclusions}
We have performed the statistical rounding error analysis for random matrix computations and provided approximate closed-form expressions for the expectation and variance of the rounding errors under a probabilistic model of the relative error. The general analysis has covered inner products, matrix-vector products, and matrix-matrix products. For specific Wishart matrices, we have presented an example of ZF detection and the corresponding LS problem and derived the approximate closed-form expressions for the expectation and variance of the rounding errors for the solution of the triangular system and LU factorization. The numerical experiments have confirmed the accuracy of the analytical expressions and demonstrated that they are generally at least two orders of magnitude tighter than alternative worst-case bounds.




\appendices
\section{Proof of {Theorem} \ref{the:inner}}
\label{app:inner}
 Assume that the sum $s_n=x_1y_1+\cdots+x_ny_n$ is evaluated from left to right \cite{higham2002accuracy}. Using the {Model} \ref{model:pro}, we have
    \begin{align*}
        \hat{s}_1&=\boldsymbol{fl}\left( x_1y_1 \right) =x_1y_1\left( 1+\delta _1 \right),\\
    \hat{s}_2&=\boldsymbol{fl}\left( \hat{s}_1+\boldsymbol{fl}\left( x_2y_2 \right) \right) =\left( \hat{s}_1+x_2y_2\left( 1+\delta _2 \right) \right) \left( 1+\delta _3 \right) \\
    &=x_1y_1\left( 1+\delta _1 \right) \left( 1+\delta _3 \right) +x_2y_2\left( 1+\delta _2 \right) \left( 1+\delta _3 \right), 
    \end{align*}
    where $\delta_i \sim dist$, $i=1,2,3$. For simplification, we denote 
    \begin{equation*}
    \begin{matrix}
         \prod_{i}^{\left( n \right)}{\left( 1+\delta _i \right)}=\underset{n\,\,\mathrm{terms}}{\underbrace{\left( 1+\delta _i \right) \cdots \left( 1+\delta _j \right) }},& i\neq j.
    \end{matrix}
    \end{equation*}
    Then, for $n = 3$, we have
    \begin{align*}
        \hat{s}_3&=\boldsymbol{fl}\left( \hat{s}_2+\boldsymbol{fl}\left( x_3y_3 \right) \right) =\left( \hat{s}_2+x_3y_3\left( 1+\delta _4 \right) \right) \left( 1+\delta _5 \right) \nonumber \\
        &=x_1y_1\prod_i^{\left( 3 \right)}{\left( 1+\delta _i \right)}+x_2y_2\prod_i^{\left( 3 \right)}{\left( 1+\delta _i \right)}+x_3y_3\prod_i^{\left( 2 \right)}{\left( 1+\delta _i \right)}.
    \end{align*}
    The pattern is clear. Overall, we have
    \begin{equation*}
        \hat{s}=\hat{s}_n=x_1y_1\prod_i^{\left( n \right)}{\left( 1+\delta _i \right)}+\sum_{k=2}^n{x_ky_k \prod_i^{\left( n-k+2 \right)}{\left( 1+\delta _i \right)} }.
    \end{equation*}
    Therefore, the rounding error $\Delta s$ is given by
    \begin{align}
        \Delta s &= \hat{s} - s \nonumber\\
        &=x_1y_1\left[ \prod_i^{\left( n \right)}{\left( 1+\delta _i \right)}-1 \right] \nonumber\\
        &+\sum_{k=2}^n{x_ky_k\left[ \prod_i^{\left( n-k+2 \right)}{\left( 1+\delta _i \right)}-1 \right]}\label{eq:delta_s}.
    \end{align}
    
    Since $x_i,y_i$ are independent of each other and $\delta_i$ have mean zero, the expectation of $\Delta s$ can be derived by using \eqref{eq:delta_s} and {Lemma} \ref{lem:e_var_product} as follows:
    {\begin{align*}
        \mathbb{E} \left( \Delta s \right) &=\mathbb{E} \left( x_1 \right) \mathbb{E} \left( y_1 \right) \mathbb{E} \left( \left[ \prod_i^{\left( n \right)}{\left( 1+\delta _i \right)}-1 \right] \right) \\
        &+\sum_{k=2}^n{\mathbb{E} \left( x_k \right) \mathbb{E} \left( y_k \right) \mathbb{E} \left( \left[ \prod_i^{\left( n-k+2 \right)}{\left( 1+\delta _i \right)}-1 \right] \right)}\\
        &=0.
    \end{align*}}
    Further, we can derive the variance of $\Delta s$ in \eqref{eq:derive_va} at the bottom of the next page, where 
    \begin{align*}
        \tau &= \mathbb{V}\left( x_ky_k \right) +\left[ \mathbb{E} \left( x_ky_k \right) \right] ^2\\
        &=\sigma _{x}^{2}\sigma _{y}^{2}+\sigma _{x}^{2}\mu _{y}^{2}+\sigma _{y}^{2}\mu _{x}^{2}+\mu _{x}^{2}\mu _{y}^{2},~ 1\leq k\leq n.
    \end{align*}
    \begin{figure*}[hb] 
    \centering 
    \hrulefill
    \begin{equation}\label{eq:derive_va}
        \begin{aligned}
        \mathbb{V}\left( \Delta s \right) &=\left\{ \mathbb{V}\left( x_1y_1 \right) +\left[ \mathbb{E} \left( x_1y_1 \right) \right] ^2 \right\} \mathbb{V}\left( \left[ \prod_i^{\left( n \right)}{\left( 1+\delta _i \right)}-1 \right] \right) 
        +\sum_{k=2}^n{\left\{ \mathbb{V}\left( x_ky_k \right) +\left[ \mathbb{E} \left( x_ky_k \right) \right] ^2 \right\} \mathbb{V}\left( \left[ \prod_i^{\left( n-k+2 \right)}{\left( 1+\delta _i \right)}-1 \right] \right)}\\
        &+2\sum_{k=2}^{n}{\mathbb{C}\left( x_1y_1\left[ \prod_i^{\left( n \right)}{\left( 1+\delta _i \right)}-1 \right] ,x_ky_k\left[ \prod_i^{\left( n-k+2 \right)}{\left( 1+\delta _i \right)}-1 \right] \right)}\\
        &+2\sum_{j=2}^{n}{\sum_{k=j+1}^{n}{\mathbb{C}\left( x_jy_j\left[ \prod_i^{\left( n-j+2 \right)}{\left( 1+\delta _i \right)}-1 \right] ,x_ky_k\left[ \prod_i^{\left( n-k+2 \right)}{\left( 1+\delta _i \right)}-1 \right] \right)}}\\
        &=\tau \left[ \left( 1+\sigma ^2 \right) ^n-1 \right] +\tau \sum_{k=2}^n{\left[ \left( 1+\sigma ^2 \right) ^{n-k+2}-1 \right]}
        +2\sum_{k=2}^n{\mathbb{E} \left( x_1y_1x_ky_k\left[ \prod_i^{\left( n-k+1 \right)}{\left( 1+\delta _{i}^{2} \right)}-1 \right] \right)}\\
        &+2\sum_{j=2}^n{\sum_{k=j+1}^n{\mathbb{E} \left( x_jy_jx_ky_k\left[ \prod_i^{\left( n-k+1 \right)}{\left( 1+\delta _{i}^{2} \right)}-1 \right] \right)}}\\
        &=\tau \left\{ \left[ \left( 1+\sigma ^2 \right) ^n-1 \right] +\sum_{k=2}^n{\left[ \left( 1+\sigma ^2 \right) ^{n-k+2}-1 \right]} \right\} 
        +2\mu _{x}^{2}\mu _{y}^{2}\sum_{i=1}^{n-1}{\left[ \frac{\left( 1+\sigma ^2 \right) \left[ \left( 1+\sigma ^2 \right) ^{n-i}-1 \right]}{\sigma ^2}-n+i \right]}\\
        &=\tau \left[ \left( 1+\sigma ^2 \right) ^n+\frac{\left( 1+\sigma ^2 \right) ^2\left[ \left( 1+\sigma ^2 \right) ^{n-1}-1 \right]}{\sigma ^2}-n \right] \\
        &+2\mu _{x}^{2}\mu _{y}^{2}\left[ \frac{\left( 1+\sigma ^2 \right) ^2\left[ \left( 1+\sigma ^2 \right) ^{n-1}-1 \right]}{\sigma ^4}-\frac{\left( n-1 \right) \left( 1+\sigma ^2 \right)}{\sigma ^2}-\frac{n\left( n-1 \right)}{2} \right].
    \end{aligned}
    \end{equation}
    \end{figure*}

    Finally,  Using Taylor expansion for $\left( 1+\sigma ^2 \right) ^n$, we have
    \begin{equation}\label{eq:ta}
        \left( 1+\sigma ^2 \right) ^n\approx 1+n\sigma ^2+\frac{n\left( n-1 \right) \sigma ^4}{2}+\frac{n\left( n-1 \right) \left( n-2 \right) \sigma ^6}{6}.
    \end{equation}
    Then, substituting \eqref{eq:ta} in \eqref{eq:derive_va}, neglecting the high-order term of $\sigma^2$ and preserving the high-order terms of $n$, we can obtain \eqref{eq:v_sigma}. And the proof ends.

\section{Proof of Corollary \ref{co:pb}}
\label{app:co}
    Note that $\Delta s = \hat{s}-s$. Using Theorem \ref{the:inner} and Lemma \ref{lem:BC}, we have 
    \begin{align}
        \left| {\Delta s}-\mathbb{E}(\Delta s) \right| &= \left| \hat{s}-s \right| \leq \sqrt{\frac{\mathbb{V} \left( \Delta s \right)}{\eta}},\label{eq:inbc}
    \end{align}
    with probability at least $1-\eta$ where $\eta = \frac{1}{\alpha^2}$ in Lemma \ref{lem:BC}. Substituting \eqref{eq:v_sigma} into \eqref{eq:inbc}, then
    \begin{align}
    \left| \hat{s}-s \right|\lesssim \sqrt{\frac{\frac{\tau}{2}n^2\sigma ^2+\frac{\mu _{x}^{2}\mu _{y}^{2}}{3}n^3\sigma ^2}{\eta}}.
    \end{align}
    
    For $\mu_x\neq0$ and $\mu_y\neq0$, we have
    \begin{align}
    \left| \hat{s}-s \right|&\overset{\left( a \right)}{\lesssim}\frac{n\left| \mu _x \right|\left| \mu _y \right|}{\sqrt{3\lambda}}\sqrt{n}\sigma =\frac{n\left| \mu _x \right|\left| \mu _y \right|}{3\sqrt{2\eta}}\sqrt{n}u \nonumber\\
    &\le \frac{\sum_{i=1}^n{\left|x_i\right|\left|y_i\right|}}{3\sqrt{2\eta}}\sqrt{n}u\nonumber\\
    &= \frac{\left| \mathbf{x} \right|^T\left| \mathbf{y} \right|}{\sqrt{\eta}}\mathcal{O} \left( \sqrt{n}u \right) \label{eq:xxx},
    \end{align}
    where $\mathcal{O}\left(n^3\sigma ^2\right)>\mathcal{O}\left(n^2\sigma ^2\right)$ at (a) point. Then \eqref{eq:bwd_n0} is obtained.

    For $\mu_x=0$ or $\mu_y=0$, without loss of generality, we consider $\mu_x=\mu_y=0$, then 
    \begin{align}
        \left| \hat{s}-s \right|&\lesssim \sqrt{\frac{\frac{\sigma _{x}^{2}\sigma _{y}^{2}}{2}n^2\sigma ^2}{\eta}}
        =\frac{n\sigma _x\sigma _y}{2\sqrt{3\eta}}u\label{eq:mean0}.
    \end{align}

    Using \cite[Corollary 3.1]{doi:10.1137/20M1314355}, if $\left|x_i\right|,\left|y_i\right|$ has bounds with means $\mu_{\left| {x} \right|}$ and $\mu_{\left| {y} \right|}$, we can obtain
    \begin{align}\label{eq:ineq}
        \left| \mathbf{x} \right|^T\left| \mathbf{y} \right|\geq \beta \mu_{\left| {x} \right|}\mu_{\left| {y} \right|}n,
    \end{align}
    where $\beta$ is a constant value (Please see the references for the details). 

    Using \eqref{eq:mean0} and \eqref{eq:ineq} yield
    \begin{align}
        \frac{\left| \hat{s}-s \right|}{\left| \mathbf{x} \right|^T\left| \mathbf{y} \right|}\lesssim \frac{\sigma _x\sigma _y}{\beta \mu _{\left| {x} \right|}\mu _{\left| {y} \right|}2\sqrt{3\eta}}u = \frac{1}{\sqrt{\eta}}\mathcal{O} \left( u \right).
    \end{align}
    
    Then the proof holds.
    
\section{Proof of {Theorem} \ref{the:tri}}
\label{app:tri}
From \eqref{eq:tri_sub}, we have
    \begin{align}          
        x_i &=\frac{b_i}{t_{ii}} -\sum_{j=1}^{i-1}\frac{t_{ij}x_j}{t_{ii}}\nonumber\\
        &\triangleq z_i - \sum_{j=1}^{i-1}c_jx_j~i=1:n.\label{eq:tri_sub_1}
    \end{align}
    Note that $t_{ii}^2\sim \chi_{m-i+1}^2,1\leq i\leq n$ and $b_i \sim \mathcal{N}\left(0,1\right)$, and we have
    \begin{align*}
        z_i&=\frac{b_i}{t_{ii}}=\frac{1}{\sqrt{m-i+1}}t,~ t\sim \mathcal{T}_{m-i+1},~1\leq i\leq n,\\
        c_j&=\frac{t_{ij}}{t_{ii}}=\frac{1}{\sqrt{m-i+1}}t,~t\sim \mathcal{T} _{m-i+1},~ 1\leq j\leq i-1,     
    \end{align*}
    
    Then, based on {Lemma} \ref{lem:e_var_student}, the expectation, variance, and covariance of $z_i~1\leq i\leq n$ and $c_j~1\leq j\leq i-1$ are given by
    \begin{align}
        \mathbb{E} \left( z_i \right) &=\mathbb{E} \left( c_j \right) =0, \label{eq:z_c_e}\\
        \mathbb{V} \left( z_i \right) &=\mathbb{V}\left( c_j \right) = \frac{1}{m-i-1}, \label{eq:z_c_var}\\
        \mathbb{C}\left( z_i,c_j \right) &=\mathbb{E} \left( \frac{b_it_{ij}}{t_{ii}^{2}} \right) -\mathbb{E} \left( z_i \right)\mathbb{E} \left( c_j \right) =0.\label{eq:z_c_cov}
    \end{align}
    Substituting \eqref{eq:z_c_e} and \eqref{eq:z_c_var} into \eqref{eq:tri_sub_1} yields
    \begin{align}
        \mathbb{E} \left( x_i \right) & = \mathbb{E} \left( z_i \right) -\sum_{j=1}^{i-1}\mathbb{E} \left( c_jx_j \right) = 0,\\
        \mathbb{V}\left( x_i \right) &=\mathbb{V}\left( z_i \right) +\sum_{j=1}^{i-1}{\mathbb{V}\left( c_j x_j \right)} \nonumber\\
        &=\frac{1}{m-i-1}+\frac{1}{m-i-1}\sum_{j=1}^{i-1}{\mathbb{V}\left( x_j \right)}\label{eq:x_var}.
    \end{align}

    Further, when \eqref{eq:tri_sub_1} is computed in floating-point arithmetic from left to right under {Model} \ref{model:pro}, we have
    \begin{align}
        \hat{x}_i&=\boldsymbol{fl}\left( \frac{\boldsymbol{fl}\left( b_i-\boldsymbol{fl}\left( \sum_{j=1}^{i-1}{\boldsymbol{fl}\left( t_{ij}\hat{x}_j \right)} \right) \right)}{t_{ii}} \right) \nonumber \\
        &=\frac{b_i\prod_k^{\left( i-1 \right)}{\left( 1+\delta _k \right)}-\sum_{j=1}^{i-1}{t_{ij}\hat{x}_j\prod_k^{\left( i-j+1 \right)}{\left( 1+\delta _k \right)}}}{t_{ii}}\left( 1+\delta _n \right) \nonumber \\
        &=\frac{b_i\prod_k^{\left( i \right)}{\left( 1+\delta _k \right)}-\sum_{j=1}^{i-1}{t_{ij}\hat{x}_j\prod_k^{\left( i-j+2 \right)}{\left( 1+\delta _k \right)}}}{t_{ii}}\label{eq:x_temp}.
    \end{align}
    Note that \eqref{eq:x_temp} involves $\hat{x}_j$ for $1\leq j\leq i-1$. To simplify, we define the relative error of $x_j$ as $\psi _j$, which has a mean of zero and variance $\sigma _{\psi _j}^{2}$, i.e.,
    \begin{equation}\label{eq:relative_x_j}
        \hat{x}_j = x_j + \Delta x_j\triangleq x_j\left( 1+\psi _j \right),~1\leq j\leq i-1.
    \end{equation}
    Utilizing {Lemma} \ref{lem:e_var_product}, we obtain
    \begin{equation}\label{eq:psi}
        \mathbb{V}\left( \Delta x_j \right) = \mathbb{V}\left( x_j \psi _j\right)=\mathbb{V}\left( x_j \right)\sigma _{\psi _j}^{2},~1\leq j\leq i-1.
    \end{equation}
    
    Then, we substitute \eqref{eq:relative_x_j} into \eqref{eq:x_temp} and have
    \begin{equation*}
        \hat{x}_i = \frac{b_i\prod_k^{\left( i \right)}{\left( 1+\delta _k \right)}-\sum_{j=1}^{i-1}{t_{ij}{x}_j\left( 1+\psi _j \right)\prod_k^{\left( i-j+2 \right)}{\left( 1+\delta _k \right)}}}{t_{ii}}.
    \end{equation*}
    Hence, the rounding error $\Delta x_i$ can be expressed as 
    \begin{align}
        \Delta x_i &= \hat{x}_i - x_i \nonumber \\
        & =\frac{b_i}{t_{ii}}\left[ \prod_k^{\left( i \right)}{\left( 1+\delta _k \right)}-1 \right] \nonumber\\
        &-\sum_{j=1}^{i-1}{\frac{t_{ij}x_j}{t_{ii}}\left[ \left( 1+\sigma _{\psi _j}^{2} \right) \prod_k^{\left( i-j+2 \right)}{\left( 1+\delta _k \right)}-1 \right]} \nonumber\\
        & =z_i\left[ \prod_k^{\left( i \right)}{\left( 1+\delta _k \right)}-1 \right]\nonumber \\
        &-\sum_{j=1}^{i-1}{c_jx_j}\left[ \left( 1+\sigma _{\psi _j}^{2} \right) \prod_k^{\left( i-j+2 \right)}{\left( 1+\delta _k \right)}-1 \right] \label{eq:delta_x}.
    \end{align}
        \begin{figure*}[hb] 
    \centering 
    \hrulefill
    \begin{equation}\label{eq:derive_tri}
        \begin{aligned}
        \mathbb{V}\left( \Delta x_i \right) &=\mathbb{V}\left( z_i \right) \mathbb{V}\left( \left[ \prod_k^{\left( i \right)}{\left( 1+\delta _k \right)}-1 \right] \right) +\sum_{j=1}^{i-1}{\mathbb{V}\left( c_j \right) \mathbb{V}\left( x_j \right) \mathbb{V}\left( \left[ \left( 1+\sigma _{\psi _j}^{2} \right) \prod_k^{\left( i-j+2 \right)}{\left( 1+\delta _k \right)}-1 \right] \right)}\\
        &\approx\frac{\left( 1+\sigma ^2 \right) ^i-1}{m-i-1}+\frac{\sum_{j=1}^{i-1}{\mathbb{V}\left( x_j \right) \left[ \left( 1+\sigma _{\psi _j}^{2} \right) \left( 1+\sigma ^2 \right) ^{i-j+2}-1 \right]}}{m-i-1}\\
        &=\frac{\left( 1+\sigma ^2 \right) ^i+\sum_{j=1}^{i-1}{\mathbb{V}\left( x_j \right) \left( 1+\sigma _{\psi _j}^{2} \right) \left( 1+\sigma ^2 \right) ^{i-j+2}}-\left( 1+\sum_{j=1}^{i-1}{\mathbb{V}\left( x_j \right)} \right)}{m-i-1}\\
        &=\frac{\left( 1+\sigma ^2 \right) ^i+\sum_{j=1}^{i-1}{\mathbb{V}\left( x_j \right) \left( 1+\sigma _{\psi _j}^{2} \right) \left( 1+\sigma ^2 \right) ^{i-j+2}}-\left( m-i-1 \right) \mathbb{V}\left( x_i \right)}{m-i-1}\\
        &=\frac{\left( 1+\sigma ^2 \right) ^i+\sum_{j=1}^{i-1}{\mathbb{V}\left( x_j \right) \left( 1+\sigma _{\psi _j}^{2} \right) \left( 1+\sigma ^2 \right) ^{i-j+2}}}{m-i-1}-\mathbb{V}\left( x_i \right).
    \end{aligned}
    \end{equation}
    \end{figure*}
    Substituting \eqref{eq:z_c_e}, \eqref{eq:z_c_var} and \eqref{eq:z_c_cov} into \eqref{eq:delta_x}, we can obtain the expectation of the rounding error $\Delta x_i$ as follows:
    \begin{align*}
        &\mathbb{E} \left( \Delta x_i \right) =\mathbb{E} \left( z_i \right) \mathbb{E} \left( \left[ \prod_k^{\left( i \right)}{\left( 1+\delta _k \right)}-1 \right] \right) \\
        &-\sum_{j=1}^{i-1}{\mathbb{E} \left( c_j \right) \mathbb{E} \left( x_j \right) \mathbb{E} \left( \left[ \left( 1+\sigma _{\psi _j}^{2} \right) \prod_k^{\left( i-j+2 \right)}{\left( 1+\delta _k \right)}-1 \right] \right)}\\
        &=0.
    \end{align*}

    Similarly, the variance of the rounding error $\Delta x_i$ is given in \eqref{eq:derive_tri} at the bottom of the page, where $\sigma _{\psi _j}^{2}$ and $\mathbb{V}\left( x_j \right)$ are given in \eqref{eq:psi} and \eqref{eq:x_var}, respectively. And the proof is done.

\section{Proof of {Theorem} \ref{the:lu}}
\label{app:lu}
First, we derive the expectation and variance of the rounding error for $u_{kj}$. When \eqref{eq:d_lu} is evaluated in float-point arithmetic from left to right under {Model} \ref{model:pro}, the computed $\hat{u}_{kj}~1\leq k \leq n,~k\leq j\leq n$ satisfies
    \begin{align}
        \hat{u}_{kj}&=\boldsymbol{fl}\left( a_{kj}-\boldsymbol{fl}\left( \sum_{i=1}^{k-1}{\boldsymbol{fl}\left( \hat{l}_{ki}\hat{u}_{ij} \right)} \right) \right)\nonumber \\
        & = a_{kj}\prod_r^{\left( k-1 \right)}{\left( 1+\delta _r \right)}-\sum_{i=1}^{k-1}{\hat{l}_{ki}\hat{u}_{ij}\prod_r^{\left( k-i+1 \right)}{\left( 1+\delta _r \right)}}\label{eq:fp_u}.
    \end{align}
    Similar to the proof of {Theorem} \ref{the:tri}, we define the relative error of $l_{ki},1\leq i\leq k-1$ and $u_{ij},1\leq i\leq k-1$ as $\epsilon _i$ and $\eta_i$, which have a mean of zero and variance $\sigma _{\epsilon _i}^{2}$ and $\sigma _{\eta_i}^{2}$, respectively. Then we have
    \begin{align}
        \hat{l}_{ki}&=l_{ki}\left( 1+\epsilon _i \right),~ 1\leq i \leq k-1,\label{eq:epsilon}\\
        \hat{u}_{ij}&=u_{ij}\left( 1+\eta _i \right),~ 1\leq i \leq k-1\label{eq:eta}.
    \end{align}
    Further, based on {Lemma} \ref{lem:e_var_product}, the variance of $\epsilon _i$ and $\eta_i$ can be expressed as 
    \begin{align}
        \sigma _{\epsilon _i}^{2}=\frac{\mathbb{V}\left( \Delta l_{ki} \right)}{\mathbb{V}\left( l_{ki} \right)},~ 1\leq i \leq k-1,\\
        \sigma _{\eta _i}^{2}=\frac{\mathbb{V}\left( \Delta u_{ij} \right)}{\mathbb{V}\left( u_{ij} \right)},~ 1\leq i \leq k-1.
    \end{align}
    Substituting \eqref{eq:epsilon} and \eqref{eq:eta} into \eqref{eq:fp_u} yields
    \begin{align*}
       \hat{u}_{kj} &= a_{kj}\prod_r^{\left( k-1 \right)}{\left( 1+\delta _r \right)} \\
       &-\sum_{i=1}^{k-1}{l_{ki}u_{ij}\left( 1+\epsilon _i \right) \left( 1+\eta _i \right) \prod_r^{\left( k-i+1 \right)}{\left( 1+\delta _r \right)}}.
    \end{align*}
    Using the definition of {Lemma} \ref{lem:term_lu_e_var}, the rounding error $\Delta u_{kj}$ can be expressed as
    \begin{align}
        \Delta u_{kj}&=\hat{u}_{kj}-u_{kj}\nonumber\\
        &=a_{kj}\left[ \prod_r^{\left( k-1 \right)}{\left( 1+\delta _r \right)}-1 \right] \nonumber\\
        &-\sum_{i=1}^{k-1}{q_{kij}\left[ \left( 1+\epsilon _i \right) \left( 1+\eta _i \right) \prod_r^{\left( k-i+1 \right)}{\left( 1+\delta _r \right)}-1 \right]}.
    \end{align}
    Note that the distributions of $a_{kj}$ and $q_{kij}$ depend on whether $j=k$ or $j\neq k$. Therefore, if $j = k$, according to {Lemma} \ref{lem:e_var_product}, {Lemma} \ref{lem:e_var_chi-square}, {Lemma} \ref{lem:dis_l_u} and {Lemma} \ref{lem:term_lu_e_var}, the expectation of $\Delta u_{kk}$ can be expressed as
    \begin{align*}
         \mathbb{E}\left( \Delta u_{kk} \right)= 0,
    \end{align*}
    and the variance of $\Delta u_{kk}$ is derived in \eqref{eq:derive_ukk} at the bottom of the page.
    \begin{figure*}[hb] 
    \centering 
    \hrulefill
    \begin{equation}\label{eq:derive_ukk}
        \begin{aligned}
        \mathbb{V}\left( \Delta u_{kk} \right) &=\mathbb{V}\left( a_{kk}\left[ \prod_r^{\left( k-1 \right)}{\left( 1+\delta _r \right)}-1 \right] \right)+\mathbb{V}\left( \sum_{i=1}^{k-1}{q_{kij}\left[ \left( 1+\epsilon _i \right) \left( 1+\eta _i \right) \prod_r^{\left( k-i+1 \right)}{\left( 1+\delta _r \right)}-1 \right]} \right) \\
        &-2\sum_{i=1}^{k-1}{\mathbb{C}\left( a_{kk}\left[ \prod_r^{\left( k-1 \right)}{\left( 1+\delta _r \right)}-1 \right] ,q_{kij}\left[ \left( 1+\epsilon _i \right) \left( 1+\eta _i \right) \prod_r^{\left( k-i+1 \right)}{\left( 1+\delta _r \right)}-1 \right] \right)}\\
        &\approx \left( m^2+2m \right) \left[ \left( 1+\sigma ^2 \right) ^{k-1}-1 \right] +3\sum_{i=1}^{k-1}{\left[ \left( 1+\sigma _{\epsilon _1}^{2} \right) \left( 1+\sigma _{\eta _1}^{2} \right) \left( 1+\sigma ^2 \right) ^{k-i+1}-1 \right]}\\
        &-2\sum_{i=1}^{k-1}{\mathbb{E} \left( a_{kk}q_{kij}\left[ \prod_r^{\left( k-i \right)}{\left( 1+\delta _{r}^{2} \right)}-1 \right] \right)}\\
        &=\left( m^2+2m \right) \left[ \left( 1+\sigma ^2 \right) ^{k-1}-1 \right] +3\sum_{i=1}^{k-1}{\left[ \left( 1+\sigma _{\epsilon _1}^{2} \right) \left( 1+\sigma _{\eta _1}^{2} \right) \left( 1+\sigma ^2 \right) ^{k-i+1}-1 \right]}
        -2\left( m+2 \right) \sum_{i=1}^{k-1}{\left[ \left( 1+\sigma ^2 \right) ^{k-i}-1 \right]}\\
        & = \left( m^2-4 \right) \left[ \left( 1+\sigma ^2 \right) ^{k-1}-1 \right] +3\left[ \sum_{i=1}^{k-1}{\left( 1+\sigma _{\epsilon _i}^{2} \right) \left( 1+\sigma _{\eta _i}^{2} \right) \left( 1+\sigma ^2 \right) ^{k-i+1}}-k+1 \right] \\
        &-2\left( m+2 \right) \left[ \frac{\left( 1+\sigma ^2 \right) \left[ \left( 1+\sigma ^2 \right) ^{k-2}-1 \right]}{\sigma ^2}-k+2 \right] .
    \end{aligned}
    \end{equation}

    \end{figure*}

    Similarly, if $j\neq k$, the expectation of $\Delta u_{kj}$ is given by 
    \begin{align*}
        &\mathbb{E}\left( \Delta u_{kj} \right)= 0,
    \end{align*}
    and the expectation of $\Delta u_{kj}$ is determined in \eqref{eq:derive_ukj} at the top the next page.
    \begin{figure*}[ht] 
    \centering 
    \begin{equation}\label{eq:derive_ukj}
        \begin{aligned}
        \mathbb{V}\left( \Delta u_{kj} \right) &=\mathbb{V}\left( a_{kj}\left[ \prod_r^{\left( k-1 \right)}{\left( 1+\delta _r \right)}-1 \right] \right) +\mathbb{V}\left( \sum_{i=1}^{k-1}{q_{kij}\left[ \left( 1+\epsilon _i \right) \left( 1+\eta _i \right) \prod_r^{\left( k-i+1 \right)}{\left( 1+\delta _r \right)}-1 \right]} \right) \\
        &-2\sum_{i=1}^{k-1}{\mathbb{C}\left( a_{kj}\left[ \prod_r^{\left( k-1 \right)}{\left( 1+\delta _r \right)}-1 \right] ,q_{kij}\left[ \left( 1+\epsilon _i \right) \left( 1+\eta _i \right) \prod_r^{\left( k-i+1 \right)}{\left( 1+\delta _r \right)}-1 \right] \right)}\\
        & \approx m\left[ \left( 1+\sigma ^2 \right) ^{k-1}-1 \right] +\sum_{i=1}^{k-1}{\left[ \left( 1+\sigma _{\epsilon _i}^{2} \right) \left( 1+\sigma _{\eta _i}^{2} \right) \left( 1+\sigma ^2 \right) ^{k-i+1}-1 \right]}
        -2\sum_{i=1}^{k-1}{\left[ \left( 1+\sigma ^2 \right) ^{k-i}-1 \right]}\\
        &=\left( m-2 \right) \left[ \left( 1+\sigma ^2 \right) ^{k-1}-1 \right] +\sum_{i=1}^{k-1}{\left( 1+\sigma _{\epsilon _i}^{2} \right) \left( 1+\sigma _{\eta _i}^{2} \right) \left( 1+\sigma ^2 \right) ^{k-i+1}}
        -2\frac{\left( 1+\sigma ^2 \right) \left[ \left( 1+\sigma ^2 \right) ^{k-2}-1 \right]}{\sigma ^2}+k-3.
    \end{aligned}
    \end{equation}
            \hrulefill
    \end{figure*}

    Then, considering the rounding error of $\bf U$, the computed $\hat{l}_{ik},1\leq i\leq k-1$ can be expressed in \eqref{eq:derive_lik} at the top of the next page, where we use \eqref{eq:epsilon}, and define $\eta_{k} = \frac{\hat{u}_{kk}-u_{kk}}{u_{kk}}$ with the mean zero and variance $\sigma _{\eta _k}^{2}=\frac{\mathbb{V}\left( \Delta u_{kk} \right)}{\mathbb{V}\left( u_{kk} \right) +\left[ \mathbb{E} \left( u_{kk} \right) \right] ^2}$ at $\left( b \right)$, and $\left( c \right)$ follows the Taylor expansion. 
    
    Further, utilizing the definition of {Lemma} \ref{lem:term_lu_e_var} and \eqref{eq:derive_lik}, the rounding error $\Delta l_{ik}$ can be expressed as
    \begin{align*}
        \Delta l_{ik}&\approx o_{ik}\left[ \left( 1-\eta _k \right) \prod_r^{\left( k \right)}{\left( 1+\delta _r \right) -1} \right] \\
        &-\sum_{j=1}^{k-1}{p_{ijk}\left[ \left( 1-\eta _k \right) \left( 1+\epsilon _j \right) \left( 1+\eta _j \right) \prod_r^{\left( k-j+2 \right)}{\left( 1+\delta _r \right)}-1 \right]}.
    \end{align*}
    \begin{figure*}[ht] 
    \centering 
    \begin{equation}\label{eq:derive_lik}
        \begin{aligned}
        \hat{l}_{ik}&=\boldsymbol{fl}\left( \frac{\boldsymbol{fl}\left( a_{ik}-\boldsymbol{fl}\left( \sum_{j=1}^{k-1}{\boldsymbol{fl}\left( \hat{l}_{ij}\hat{u}_{jk} \right)} \right) \right)}{\hat{u}_{kk}} \right),~k+1\leq i \leq n \\
        &=\frac{a_{ik}\prod_r^{\left( k-1 \right)}{\left( 1+\delta _r \right)}-\sum_{j=1}^{k-1}{\hat{l}_{ij}\hat{u}_{jk}\prod_r^{\left( k-j+1 \right)}{\left( 1+\delta _r \right)}}}{\hat{u}_{kk}}\left( 1+\delta _o \right) \\
        &\overset{\left( b \right)}{=}\frac{a_{ik}\prod_r^{\left( k \right)}{\left( 1+\delta _r \right)}-\sum_{j=1}^{k-1}{l_{ij}u_{jk}\left( 1+\sigma _{\epsilon _j}^{2} \right) \left( 1+\sigma _{\eta _j}^{2} \right) \prod_r^{\left( k-j+2 \right)}{\left( 1+\delta _r \right)}}}{u_{kk}\left( 1+\eta _k \right)}\\
        &\overset{\left( c \right)}{\approx}\frac{a_{ik}\prod_r^{\left( k \right)}{\left( 1+\delta _r \right)}-\sum_{j=1}^{k-1}{l_{ij}u_{jk}\left( 1+\sigma _{\epsilon _j}^{2} \right) \left( 1+\sigma _{\eta _j}^{2} \right) \prod_r^{\left( k-j+2 \right)}{\left( 1+\delta _r \right)}}}{u_{kk}}\left( 1-\eta _k \right). 
    \end{aligned}
    \end{equation}
    \hrulefill
    \end{figure*}

    \begin{figure*}[ht] 
    \centering 
    {\begin{equation}\label{eq:derive_lik_v}
        \begin{aligned}
        \mathbb{V}\left( \Delta l_{ik} \right) &\approx \mathbb{V}\left( o_{ik}\left[ \left( 1-\eta _k \right) \prod_r^{\left( k \right)}{\left( 1+\delta _r \right) -1} \right] \right) 
        +\mathbb{V}\left( \sum_{j=1}^{k-1}{p_{ijk}\left[ \left( 1-\eta _k \right) \left( 1+\epsilon _j \right) \left( 1+\eta _j \right) \prod_r^{\left( k-j+2 \right)}{\left( 1+\delta _r \right)}-1 \right]} \right) \\
        &-2\sum_{j=1}^{k-1}{\mathbb{C}\left( o_{ik}\left[ \left( 1-\eta _k \right) \prod_r^{\left( k \right)}{\left( 1+\delta _r \right) -1} \right] ,p_{ijk}\left[ \left( 1-\eta _k \right) \left( 1+\epsilon _j \right) \left( 1+\eta _j \right) \prod_r^{\left( k-j+2 \right)}{\left( 1+\delta _r \right)}-1 \right] \right)}\\
        &\approx \frac{\left(m-4\right)\left[ \left( 1+\sigma _{\eta _k}^{2} \right) \left( 1+\sigma ^2 \right) ^k-1 \right]}{\left( m-k-1 \right) \left( m-k-3 \right)} +\frac{\sum_{j=1}^{k-1}{\left[ \left( 1+\sigma _{\eta _k}^{2} \right) \left( 1+\sigma _{\epsilon _j}^{2} \right) \left( 1+\sigma _{\eta _j}^{2} \right) \left( 1+\sigma ^2 \right) ^{k-j+2}-1 \right]}}{\left( m-k-1 \right) \left( m-k-3 \right)}\\
        &-2\sum_{j=1}^{k-1}{\mathbb{E} \left( o_{ik}p_{ijk}\left[ \left( 1+\eta _{k}^{2} \right) \prod_i^{\left( k-j+1 \right)}{\left( 1+\delta _{i}^{2} \right) -1} \right] \right)}\\
        &=\frac{\left(m-4\right)\left[ \left( 1+\sigma _{\eta _k}^{2} \right) \left( 1+\sigma ^2 \right) ^k-1 \right]}{\left( m-k-1 \right) \left( m-k-3 \right)} +\frac{\sum_{j=1}^{k-1}{\left[ \left( 1+\sigma _{\eta _k}^{2} \right) \left( 1+\sigma _{\epsilon _j}^{2} \right) \left( 1+\sigma _{\eta _j}^{2} \right) \left( 1+\sigma ^2 \right) ^{k-j+2}-1 \right]}}{\left( m-k-1 \right) \left( m-k-3 \right)}\\
        &-\frac{2\sum_{j=1}^{k-1}{\left[ \left( 1+\sigma _{\eta _k}^{2} \right) \left( 1+\sigma ^2 \right) ^{k-j+1}-1 \right]}}{\left( m-k-1 \right) \left( m-k-3 \right)}\\
        &=\frac{\left(m-6\right)\left[ \left( 1+\sigma _{\eta _k}^{2} \right) \left( 1+\sigma ^2 \right) ^k-1 \right]}{\left( m-k-1 \right) \left( m-k-3 \right)}
        +\frac{\left( 1+\sigma _{\eta _k}^{2} \right) \sum_{i=j}^{k-1}{\left( 1+\sigma _{\epsilon _j}^{2} \right) \left( 1+\sigma _{\eta _j}^{2} \right) \left( 1+\sigma ^2 \right) ^{k-j+2}}-k+1}{\left( m-k-1 \right) \left( m-k-3 \right)}\\
        &-\frac{2}{\left( m-k-1 \right) \left( m-k-3 \right)}\left[ \left( 1+\sigma _{\eta _k}^{2} \right) \frac{\left( 1+\sigma ^2 \right) ^2\left[ \left( 1+\sigma ^2 \right) ^{k-2}-1 \right]}{\sigma ^2}-k+2 \right].  
    \end{aligned}
    \end{equation}}
    \hrulefill
    \end{figure*}

    Therefore, based on {Lemma} \ref{lem:e_var_product}, {Lemma} \ref{lem:e_var_wishart}, {Lemma} \ref{lem:dis_l_u} and {Lemma} \ref{lem:term_lu_e_var}, the variance of $\Delta l_{ik}$ can be derived in \eqref{eq:derive_lik_v} at the top of the next page. Therefore, the theorem holds.

\bibliographystyle{IEEEtran}
\bibliography{reference}

\end{document}